%% file: GBundlesProjectiveLine.tex
\documentclass[a4paper,11pt, leqno]{article}
\usepackage[T1]{fontenc}
\usepackage{articlemacro}

\newcommand{\B}[1]{\mathbb{B}#1}
\newcommand{\Bunline}{\underline{\mathcal{B}\textup{un}}}
\DeclareMathOperator{\ClOpen}{OC}

\begin{document}

\title{Extension and lifting of $G$-bundles for stacks}

\author{Torsten Wedhorn}

\maketitle


\noindent{\scshape Abstract.\ } We study extension properties for morphisms of stacks of bundles for group algebraic spaces. Applications are a short proof of the classification of bundles on the projective line for affine smooth group schemes whose identity component is reductive and the existence of splittings of filtered fiber functors for affine smooth group schemes defined over affine schemes.

\noindent{\scshape MSC.\ } 14D23 (Primary), 14L15, 14L30, 14N05 (Secondary)


\input{Ch0GBundlesProjectiveLine.tex}

\input{Ch1GBundlesProjectiveLine.tex}

\input{Ch2GBundlesProjectiveLine.tex}

\input{Ch3GBundlesProjectiveLine.tex}

%


\bibliographystyle{plain}
\bibliography{references}

\noindent \textit{Torsten Wedhorn}: \textsc{TU Darmstadt, Department of Mathematics, Germany}

\noindent\textit{Email}: \texttt{wedhorn@mathematik.tu-darmstadt.de}

\end{document}

%% file: Ch0GBundlesProjectiveLine.tex

\section*{Introduction}

In this paper we will study extension properties of $G$-bundles for morphisms of (not necessarily algebraic) stacks $\Xcal$, where $G$ is a group scheme (or even a group algebraic space). We denote the groupoid of $G$-bundles on $\Xcal$ by $\Bun_G(\Xcal)$ and the set of isomorphism classes in $\Bun_G(\Xcal)$ by $H^1(\Xcal,G)$. The classifying stack of $G$ is denoted by $\B{G}$. Before giving the main results, let us start with an illustrative application.

\begin{intro-theorem}[Theorem~\ref{ClassifyGBundlesPP1}]\label{IntroGBunPP^1}
Let $k$ be a field and let $G$ be an affine smooth group scheme over $k$ such that its identity component is reductive. Then there is a natural bijection
\[
H^1(\PP^1_k,G) \cong H^1(\B{\GG_{m,k}},G),
\]
where $\B{\GG_{m,k}}$ denotes the classifying stack of the multiplicative group over $k$.
\end{intro-theorem}

The right hand side can be identified (Section~\ref{SecGBundleBGGm}) with isomorphism classes of pairs $(P,\mu)$, where $P$ is a $G$-bundle over $k$ and $\mu\colon \GG_{m,k} \to \Autline_G(P)$ is a cocharacter of the strong inner form of $G$ attached to $P$. The classification of $G$-bundles over $\PP^1_k$ has a long history. For reductive groups $G$, Grothendieck \cite{Groth_ClassFibresHol} gave a classification over an algebraically closed field which was generalized by Harder \cite[Satz~3.4]{Harder_HeGruppenCurves} to a classification of Zariski locally trivial $G$-bundles over any field for a split reductive group $G$. Biswas and Nagaraj classified in \cite{BiswasNagaraj} Zariski locally trivial $G$-bundles over $\PP^1_k$ for any field and any reductive group over $k$. Ansch\"utz proved Theorem~\ref{IntroGBunPP^1} in \cite{Anschuetz_TannakianClassPP1} for all fields $k$ and all reductive groups $G$ expressing the right hand side in a ``Tannakian way'' (Remark~\ref{BunGBGGmTannaka}).

The main interest of Theorem~\ref{IntroGBunPP^1} does not lie in its generalization to not necessarily connected groups (although it does have the amusing side effect of giving a new proof of $\pi_1(\PP^1_k) \cong \pi_1(k)$, see Remark~\ref{Pi1PP1}), but in its method of proof. In \cite{BiswasNagaraj} and \cite{Anschuetz_TannakianClassPP1} the proof was obtained by using the Harder--Narasimhan stratification and the classification of vector bundles over $\PP^1_k$. Here we bypass such arguments completely and proceed as follows. Consider the diagram of algebraic stacks
\[
\PP^1_k = [\GG_{m,k}\backslash (\AA^2_k \setminus \{0\})] \ltoover{j} [\GG_{m,k}\backslash \AA^2_k] \llefttoover{i} [\GG_{m,k}\backslash \{0\}] = \B{\GG_{m,k}}.
\]
Via pullback we obtain maps
\[
H^1(\PP^1_k,G) \llefttoover{j^*} H^1([\GG_{m,k}\backslash \AA^2_k],G) \ltoover{i^*} H^1(\B{\GG_{m,k}},G)\tag{$\spadesuit$}
\]
and as immediate applications of our main results we obtain that $j^*$ and $i^*$ are bijections.

The main focus of the paper is the following question. For a morphism $f\colon \Ycal \to \Xcal$ of stacks (not necessarily algebraic) we study properties of the pullback functor $f^*\colon \Bun_G(\Xcal) \to \Bun_G(\Ycal)$ and aim to give criteria when
\begin{assertionlist}
\item[(I)]
the functor $f^*$ is fully faithful or even an equivalence of categories, or
\item[(II)]
$f$ is a closed immersion and $f^*$ induces a surjection or bijection on isomorphism classes.
\end{assertionlist}
To be more concrete, let us introduce some notation. Let $S$ be an algebraic space and let $G$ be a group algebraic space over $S$. In the introduction, all stacks are defined over $S$, are stacks for the fpqc-topology, have a diagonal representable by algebraic spaces, and admit an fpqc-covering by an algebraic space. Regarding Problem~(I), we prove the following general result.

\begin{intro-proposition}[Proposition~\ref{RestrictionFullyFaithful}, Proposition~\ref{InEssentialImage}]\label{IntroExtendGBundle}
Let $f\colon \Ucal \to \Xcal$ be a qcqs morphism of stacks that is representable by schemes. Suppose that $\Oscr_{\Xcal} = f_*\Oscr_{\Ucal}$ (a condition that can be checked fpqc locally on $\Xcal$).
\begin{assertionlist}
\item
If $G \to S$ is affine and flat, then $f^*\colon \Bun_G(\Xcal) \to \Bun_G(\Ucal)$ is fully faithful.
\item
Suppose that there exists an embedding $\iota\colon G \mono H$ of flat affine group algebraic spaces such that the fppf-quotient $G\backslash H$ is representable by an algebraic space that is affine over $S$. Then a $G$-bundle $P$ on $\Ucal$ is in the essential image of $f^*$ if and only if the $H$-bundle $\iota_*(P)$ is fpqc locally on $\Ucal$ in the essential image of $f^*\colon \Bun_H(\Xcal) \to \Bun_H(\Ucal)$.
\end{assertionlist}
\end{intro-proposition}

Such questions have also been considered in many other papers in various special situations (e.g., \cite[Th\'eor\`eme~6.13]{ColTheSans_FibreQuadratiques}, \cite[Section 3.5]{Anschuetz_ExtendingTorsors} to name just two). In this generality, the result seems to be new, even though it is certainly not surprising.

We also give applications of Proposition~\ref{IntroExtendGBundle} if $f$ is an open immersion (Proposition~\ref{ExtendGBundleOpenImmersion}), or if $f$ is proper, flat, of finite presentation with geometrically connected and geometrically reduced fibers (Proposition~\ref{PullBackProper}). In each case we choose $H$ a general linear group and use known results for vector bundles. For instance, if $f$ is the open embedding $j\colon \PP_k^1 \to [\GG_m\backslash \AA_k^2]$, then fppf-locally $f$ is the embedding $U := \AA^2_k \setminus \{0\} \to \AA^2_k$ and by the Auslander--Buchsbaum theorem, every vector bundle on $U$ can be extended uniquely to $\AA^2_k$ since $\AA^2_k$ is regular of dimension $2$. Hence Theorem~\ref{IntroExtendGBundle} immediately implies that the map $j^*$ in ($\spadesuit$) is bijective.

To approach Problem~(II), we consider pairs consisting of an algebraic stack and a closed substack that are henselian (Definition~\ref{DefHenselianPair}) and apply a result of Alper, Hall, and Rydh \cite[Proposition~7.10]{AHR_EtaleLocal} to obtain the following result.

\begin{intro-theorem}[Theorem~\ref{LiftHenselian}]\label{IntroLiftHenselian}
Let $A$ be a ring and let $H$ be a group of multiplicative type of finite type over $A$ that becomes diagonalizable after a finite locally free surjective base change. Suppose that $H$ acts on an affine $A$-scheme $\Xtilde$ and let $\Ztilde$ be an $H$-equivariant closed subscheme. Then $([H\backslash \Xtilde], [H\backslash \Ztilde])$ is a henselian pair if and only if the fixed point loci $(\Xtilde^H,\Ztilde^H)$ form a henselian pair. In this case, one has
\[
H^1([H\backslash \Xtilde],G) \iso H^1([H\backslash \Ztilde],G)
\]
for every smooth quasi-affine group scheme $G$ over $A$.
\end{intro-theorem}

If $H$ is trivial, the theorem has been obtained by Bouthier and \v{C}esnavi\v{c}ius \cite[Theorem~2.1.6]{BC_Torsors}.

In fact, the theorem also holds for the groupoids of morphisms into any gerbe that is fpqc-locally the classifying stack of a smooth quasi-affine group scheme.

Along the way, we prove in Section~\ref{Sec:ClopenConnComp} the following purely topological result, where $\ClOpen(X)$ denotes the set of open and closed subsets of a topological space $X$.

\begin{intro-proposition}[Proposition~\ref{ConnMapPreSpektral}]\label{IntroClopenPi0}
Let $X$ and $Y$ be spectral topological spaces (e.g., the underlying topological spaces of qcqs algebraic stacks \cite[0DQP]{Stacks}). Let $f\colon X \to Y$ be a continuous map. Then $f$ induces a bijection $\ClOpen(Y) \iso \ClOpen(X)$ if and only if $\pi_0(f)\colon \pi_0(X) \to \pi_0(Y)$ is bijective.
\end{intro-proposition}

This result allows to give a new characterization of henselian pairs of algebraic stacks, see Proposition~\ref{PropHenselianMorph}.

If $H$ is diagonalizable, we give explicit criteria when $([H\backslash \Xtilde], [H\backslash \Ztilde])$ is a henselian pair in Section~\ref{SecHenselianGraded}. For instance, these criteria show that if $B$ is a $\ZZ$-graded $A$-algebra concentrated in degrees $\geq 0$, then $([\GG_m\backslash (\Spec B)],\B{\GG_{m,B_0}})$ is a henselian pair (Example~\ref{TrivialHenselian}). In particular $[\GG_{m,A}\backslash \AA_A^n, \B{\GG_{m,A}}]$ is a henselian pair and hence Theorem~\ref{IntroLiftHenselian} implies for $n = 2$ that the map $i^*$ in ($\spadesuit$) is bijective.

This strategy to prove Theorem~\ref{IntroGBunPP^1} has been further extended in \cite{BCS_TorsorsConstant} to obtain Theorem~\ref{IntroGBunPP^1} for smooth $k$-group schemes $G$ whose largest connected, smooth, affine $k$-subgroup has no nontrivial split unipotent normal $k$-subgroups, see \cite[Theorem~1.3]{BCS_TorsorsConstant}.

For $n = 1$, Theorem~\ref{IntroLiftHenselian} also yields a bijection
\[
H^1([\GG_{m,A}\backslash \AA^1_A],G) \iso H^1(\B{\GG_{m,A}},G)\tag{*}
\]
for every quasi-affine smooth group scheme $G$ over $A$. If $A$ is an $A_0$-algebra for some ring $A_0$ and $G$ is the base change of a group scheme $G_0$ over $A_0$ such that $\B{G_0}$ satisfies the resolution property (Example~\ref{BGResolution}), then $G$-bundles over a stack $\Xcal$ are given by exact tensor functors from the category $\Rep(G_0)$ of $G$-representations on finite projective $A_0$-modules to the category of vector bundles on $\Xcal$ (Proposition~\ref{TannakaGBundles}). As recalled in Section~\ref{SecSplitFiberFunctor}, we can identify the left hand side of (*) with the set of isomorphism classes of filtered fiber functors of $\Rep(G_0)$ over $A$ and the right hand side with the set of isomorphism classes of graded fiber functors of $\Rep(G_0)$ over $A$. Hence we obtain as a corollary of Theorem~\ref{IntroLiftHenselian} that every filtered fiber functor of $\Rep(G)$ admits a splitting (Corollary~\ref{SplittingFFF}).
This generalizes a theorem of Ziegler \cite{Ziegler_FFF} for neutralizable Tannakian categories from the case that $A$ is a field to a general ring. Ziegler also gave in \cite{Ziegler_FFFArbitrary} independently a different proof of the splitting of filtered fiber functors in the situation above if $G$ is in addition affine.

\bigskip

\noindent{\scshape Acknowledgements.\ }

I would like to thank Timo Richarz and Paul Ziegler for helpful comments. I am very grateful to the referee whose many remarks led to a significant improvement of the paper. In particular, I thank the referee for the remark that Theorem~\ref{LiftHenselian} should be generalized from affine to quasi-affine group schemes. This project was supported by the German Research Foundation (DFG), TRR 326.


\tableofcontents

\bigskip\bigskip

\subsection*{Notation}

All rings are commutative with $1$ unless explicitly stated otherwise.

Let $S$ be an algebraic space. We denote by $\Affrel{S}$ the full subcategory of the category of $S$-schemes consisting of morphisms $X \to S$ of algebraic spaces with $X$ an affine scheme. A (classical) prestack over $S$ is a contravariant functor of $(2,1)$-categories from the category $\Affrel{S}$ to the $(2,1)$-category of groupoids. Equivalently, a prestack is a presheaf on $\Affrel{S}$ with values in 1-truncated anima. If $\Scal$ is a prestack, then we also denote by $\Affrel{\Scal}$ the category of morphisms $X \to \Scal$ of prestacks with $X$ an affine scheme.

Let $\Xcal$ and $\Ycal$ be prestacks over a prestack $\Scal$. We denote by $\Hom_{\Scal}(\Xcal,\Ycal)$ the groupoid of morphisms $\Xcal \to \Ycal$ of prestacks over $\Scal$. The prestack $T \sends \Hom_{\Scal}(\Xcal \times_{\Scal} T,\Ycal)$ is denoted by $\Homline_{\Scal}(\Xcal,\Ycal)$. This is a stack, if $\Xcal$ and $\Ycal$ are stacks.

A morphism $f\colon \Xcal \to \Ycal$ of prestacks is called  \emph{representable by schemes} (resp.~\emph{by algebraic spaces}) if for every morphism $Y \to \Ycal$ where $Y$ is a scheme, the fiber product $\Xcal \times_{\Ycal} Y$ is a scheme (resp.~an algebraic space). If $f$ is representable by schemes (resp.~algebraic spaces) and $\Pbf$ is a property of scheme morphisms (resp.~of morphisms of algebraic spaces) that is stable under base change in the category of schemes (resp.~of algebraic spaces), then $f$ is said to have $\Pbf$ if for every morphism $Y \to \Ycal$, where $Y$ is a scheme (resp.~$Y$ is an algebraic space), the base change $\Xcal \times_{\Ycal} Y \to Y$ has property $\Pbf$. 

The notion of gerbe is always meant with respect to the fpqc topology, i.e. a morphism of algebraic stacks $f\colon \Ycal \to \Xcal$ is called a gerbe if it is fpqc-locally on $\Xcal$ isomorphic to a classifying stack for some fpqc-sheaf of groups.

If $\Xcal$ is a prestack,
we denote the category of vector bundles on $\Xcal$ by $\Vec{\Xcal}$. It is endowed with the standard exact structure. The usual tensor product makes the category of $\Oscr_{\Xcal}$-modules into a symmetric monoidal category.
It induces by restriction a symmetric monoidal structure on $\Vec{\Xcal}$.

%% file: Ch1GBundlesProjectiveLine.tex
\section{Preliminaries on quotient stacks and $G$-bundles on stacks}

\subsection{Quotient stacks for the fpqc-topology}

Let $S$ be an algebraic space. We denote by $\Affrel{S}$ the category of affine schemes endowed with a morphism to $S$. We will deal systematically with sheaves (of sets or of groupoids) for the fpqc topology on $\Affrel{S}$ (or equivalently on the category of all algebraic spaces over $S$) since we will deal with general flat group schemes (not necessarily of finite presentation)

The category $\Affrel{S}$ is not small. Moreover, for any given non-empty affine scheme $T$ over $S$ the category of fpqc-covering sieves of $T$ does not have a small cofinal subcategory. Hence categories of sheaves in the naive sense fail to be Grothendieck topoi and presheaves do not admit a sheafification in general. To overcome this difficulty we proceed by first limiting cardinalities and then forming the colimit, following \cite[Section~2.1 and App.~A]{HesselholtPstragowski_DiracGeometryII} where this is explained in the more complicated case of the fpqc topology on the category of Dirac rings, see also \cite[\S2]{Gregoric_StoneDuality}.

We will consider presheaves and sheaves with value in arbitrary anima (= spaces = $\infty$-groupoids) even though in this paper we need only (pre)sheaves with values in 1-truncated anima (= groupoids). We denote the $\infty$-category of anima by $\Ani$.

Let $\kappa$ be an uncountable strong limit cardinal and let $R$ be a $\kappa$-small ring. Let $\Aff_{R,\kappa}$ be the opposite category of $\kappa$-small $R$-algebras, which we always endow with the fpqc-topology. Let $\Sh{\Aff_{R,\kappa}}$ be the $\infty$-category of fpqc-sheaves on $\Aff_{R,\kappa}$ with values in $\Ani$, see, for instance, \cite[Proposition 2.5]{HesselholtPstragowski_DiracGeometryII}.

If $\kappa' > \kappa$ are strongly inaccessible cardinals, then left Kan extension along $\Aff_{R,\kappa} \mono \Aff_{R,\kappa'}$ induces fully faithful functors $\Sh{\Aff_{R,\kappa}} \mono \Sh{\Aff_{R,\kappa'}}$ that commute with all colimits and with $\kappa$-small limits \cite[Theorem 2.15, Lemma A.8]{HesselholtPstragowski_DiracGeometryII}. For any ring $R$ we set
\[
\Sh{\Aff_{R}} := \colim_{\kappa}\Sh{\Aff_{R,\kappa}},
\]
where the colimit is indexed by the large partially ordered set of uncountable strongly inaccessible cardinals $\kappa$ such that $R$ is $\kappa$-small.

Then $\Sh{\Aff_{R}}$ is equivalent to the $\infty$-category of fpqc-sheaves on $\Aff_R$ which are accessible as a presheaf \cite[Remark 2.16]{HesselholtPstragowski_DiracGeometryII} or, equivalently by \cite[Proposition A.2]{HesselholtPstragowski_DiracGeometryII}, which are as presheaves a colimit of a small diagram of representable functors. We denote the $\infty$-category of these presheaves by $\PSh{\Aff_{R}}$.

By \cite[Corollary 2.17]{HesselholtPstragowski_DiracGeometryII} we obtain:

\begin{proposition}\label{CritExistenceSheafification}
The inclusion $\Sh{\Aff_R} \to \PSh{\Aff_R}$ admits a left adjoint sheafification functor which commutes with finite limits. 
\end{proposition}

Now let $G$ be an fpqc-sheaf of groups on the category of $S$-schemes. For instance, $G$ could by any group algebraic space over $S$ by a theorem of Gabber \cite[0APL]{Stacks}. Let $X$ be an fpqc-sheaf over $S$ with a left $G$-action $a\colon G \times_S X \to X$. We define $[G\backslash X]$ as the fpqc-stack in groupoids on $\Ccal$ attached to the prestack
\[
[G\,{}_p\!\backslash X]\colon T \sends [G(T)\backslash X(T)]
\]
where $[G(T)\backslash X(T)]$ is the groupoid whose objects are the elements of $X(T)$ and whose maps $x \to x'$ are the elements $g \in G(T)$ with $gx = x'$. Composition is defined via the group law in $G(T)$. The stackification (i.e. sheafification of presheaves with values in 1-truncated anima) of $[G\,{}_p\!\backslash X]$ always exists by Proposition~\ref{CritExistenceSheafification} since we can write $[G\,{}_p\!\backslash X]$ as the (homotopy) colimit of the bar construction
\begin{equation}\label{EqBarResolution}
\xymatrix{\dots \ar@<0.5ex>[r] \ar@<1.5ex>[r] \ar@<-0.5ex>[r] \ar@<-1.5ex>[r] & G \times_S G \times_S X \ar@<1ex>[r] \ar[r] \ar@<-1ex>[r] & G \times_S X \ar@<0.5ex>[r]^-a \ar@<-0.5ex>[r]_-{p_2} & X.}
\end{equation}
As usual, the quotient stack can also be described via $G$-bundles, which also shows its existence. Here we mean by a $G$-bundle an fpqc-sheaf $P$ on $\Affrel{S}$ endowed with a $G$-action that is fpqc-locally isomorphic to $G$ endowed with the $G$-action by left multiplication. Then $[G\backslash X]$ is identified with the stack that sends an object $T$ of $\Affrel{S}$ to the groupoid of pairs $(P,\varphi)$, where $P$ is a $G$-bundle over $T$ and $\varphi\colon P \to X \times_S T$ is a $G$-equivariant map of sheaves over $T$. A morphism $(P,\varphi) \to (P',\varphi')$ in $[G\backslash X](T)$ is a $G$-equivariant isomorphism $u\colon P \iso P'$ of sheaves such that $\varphi = \varphi' \circ u$. 

%
%

In particular, we have the classifying stack $\B{G} := [G\backslash S]$, where $G$ acts trivially on $S$, which classifies $G$-bundles for the fpqc-topology.

In the sequel we will also use the following notion. Recall from \cite[0AP5]{Stacks} that a morphism $f\colon X \to Y$ of schemes is called \emph{ind-quasi-affine} if for every open affine subscheme $V$ of $Y$ every open quasi-compact subscheme of $f^{-1}(V)$ is quasi-affine. The property to be ind-quasi-affine is stable under composition and base change and stable under fpqc-descent. Every immersion and every quasi-affine morphism is ind-quasi-affine. By the usual argument (e.g., \cite[Remark~9.11]{GW1}) this implies that if a composition $g \circ f$ of scheme morphisms is ind-quasi-affine, then $f$ is ind-quasi-affine.

We call a morphism $\Xcal \to \Ycal$ of prestacks \emph{ind-quasi-affine} if it is representable by schemes and for every morphism of schemes $T \to \Ycal$ the base change $\Xcal \times_{\Ycal} T \to T$ is ind-quasi-affine.

\begin{remark}\label{PropertyGBundlesIndQuasiAffine}
Let $G$ be an algebraic group space over $S$ such that $G \to S$ is ind-quasi-affine and an fpqc covering.
\begin{assertionlist}
\item\label{PropertyGBundlesIndQuasiAffine1}
If $X$ is an algebaic space over $S$ endowed with a $G$-action, then the quotient map $X \to [G\backslash X]$ is also an fpqc covering and ind-quasi-affine.

More generally, the arguments in \cite[04X0]{Stacks} show that if $G \to S$ has a property $\Pbf$ that can be checked fpqc-locally on the target. Then the canonical map $X \to [G\backslash X]$ has $\Pbf$.
\item
If $P$ is an fpqc-stack on $\Affrel{S}$ which is fpqc-locally isomorphic to $G$ (for instance, a $G$-bundle for the fpqc-topology), then the structure map $P \to S$ is representable by schemes.

Indeed since $P$ is an fpqc-stack and fpqc-locally has values in sets, as this is the case for $G$, the stack $P$ itself has values in sets, i.e., $P$ is an fpqc-sheaf. Hence the claim follows since one has effective fpqc descent for ind-quasi-affine schemes by a result of Gabber \cite[0APK]{Stacks}.
\item
If $G$ is in addition flat and locally of finite presentation over $S$, then every $G$-bundle $P$ over $T$ is faithfully flat and locally of finite presentation over $T$, and we can choose $P$ itself as a covering and obtain that every $G$-bundles is already fppf-locally trivial.

If $G$ is even smooth over $S$, then every $G$-bundle is even \'etale locally trivial, since every smooth morphism has \'etale locally sections.
\end{assertionlist}
\end{remark}

\begin{proposition}\label{QuotientByIndQuasiAffine}
Let $S$ be an algebraic space, let $G$ be a group algebraic space over $S$ such that $G \to S$ is an fpqc-covering and ind-quasi-affine. Let $X$ be an algebraic space over $S$ endowed with an action of $G$. Then the diagonal $\Delta_{[G\backslash X]}$ is ind-quasi-affine, and in particular representable by schemes.

If $X \to S$ has an affine diagonal and $G \to S$ is affine, then $\Delta_{[G\backslash X]}$ is affine.
\end{proposition}

\begin{proof}
Consider the following diagram
\begin{equation}\label{EqDiagramDiagonal}
\begin{aligned}\xymatrix{
G \times_S X \ar[r]^{\sigma} \ar[d] & X \times_S X \ar[r]^{p_2} \ar[d]^p & X \\
[G\backslash X] \ar[r]^-{\Delta_{[G\backslash X]}} & [G\backslash X] \times_S [G\backslash X],
}\end{aligned}
\end{equation}
where the square is 2-cartesian. Here $\sigma$ is the map given by $(g,x) \sends (gx,x)$ on points. Since $G \to S$ is ind-quasi-affine, so is $p_2 \circ \sigma$. Therefore, $\sigma$ is ind-quasi-affine by cancellation which implies that $\Delta_{[G\backslash X]}$ is ind-quasi-affine by fpqc-descent.

For ``affine'' the argument is the same where we use that if a composition $g \circ f$ is affine and $g$ has affine diagonal, then $f$ is affine.
\end{proof}

\begin{remark}\label{DiagonalClassifying}
The 2-cartesian square in \eqref{EqDiagramDiagonal} applied to $X = S$ shows that if $G \to S$ is an fpqc-covering and $\Pbf$ is a property stable under fpqc base change and under fpqc-descent, then $\Delta_{\B{G}}$ has $\Pbf$ if and only if $G \to S$ has $\Pbf$.
\end{remark}


\subsection{$G$-bundles over stacks}

We continue to assume that $S$ is an algebraic space.

\begin{defrem}\label{DefGBundles}
Let $\Xcal$ be a prestack over $S$ and let $G$ be a group algebraic space over $S$. We define the \emph{groupoid $\Bun_G(\Xcal)$ of $G$-bundles over $\Xcal$} as the groupoid of $S$-morphisms $\Xcal \to \B{G}$, i.e., we set
\[
\Bun_G(\Xcal) := \Hom_S(\Xcal,\B{G}).
\]
We denote by $H^1(\Xcal,G)$ the set of isomorphism classes of $G$-bundles on $\Xcal$.

The canonical map $S \to \B{G}$ corresponds to the trivial $G$-bundle on $S$. More generally, the morphism $\Xcal \to S \to \B{G}$ is called the \emph{trivial $G$-bundle on $\Xcal$}. By distinguishing it, $H^1(\Xcal,G)$ is made into a pointed set.

We also denote by $\Bunline_G(\Xcal)$ the prestack over $S$ given by
\[
\Bunline_G(\Xcal)(T) = \Hom_S(\Xcal \times_S T,\B{G}).
\]
If $u\colon \Ycal \to \Xcal$ is a morphism of prestacks over $S$ and $\Escr$ is a $G$-bundle on $\Xcal$ given by a map $\Xcal \to \B{G}$, then its \emph{pullback by $u$} is the $G$-bundle on $\Ycal$ given by the composition $\Ycal \ltoover{u} \Xcal \to \B{G}$. It is denoted by $u^*\Escr$.
\end{defrem}

Suppose that $G \to S$ is ind-quasi-affine and an fpqc-covering. In many cases, it is possible to describe $G$-bundles via their representations. 
We set
\[
\Rep_S(G) := \Rep(G) := \Vec{\B{G}}
\]
which we consider as an exact symmetric monoidal category. Since $S \to \B{G}$ is an fpqc-covering and fpqc-descent is effective for vector bundles, this is the category of vector bundles $\Escr$ over $S$ endowed with a homomorphism of $S$-group algebraic spaces $G \to \Autline_{\Oscr_S}(\Escr)$, where $\Autline_{\Oscr_S}(\Escr)$ is the reductive group scheme over $S$ whose $T$-valued points for an affine $S$-scheme $f\colon T \to S$ are given by the group of $\Oscr_T$-linear automorphisms of $f^*\Escr$.

Then we have the following application of a result of Tonini \cite[Theorem~A]{Tonini_Tannaka}.

\begin{proposition}\label{TannakaGBundles}
Let $S$ be a quasi-compact scheme and let $G$ be a group scheme over $S$ such that $G \to S$ is flat and quasi-affine. Suppose that $\B{G}$ has the resolution property, i.e., every quasi-coherent module over $\B{G}$ is a quotient of a direct sum of vector bundles. Let $\Xcal$ be any prestack over $S$. Then the canonical morphism of groupoids
\begin{equation}\label{EqToniniTannaka}
\begin{aligned}
\Bun_G(\Xcal) = \Hom_S(\Xcal,\B{G}) &\lto \left\{\begin{matrix}\text{exact strong monoidal $\Oscr_S$-linear}\\ \text{functors $\Rep(G) \to \Vec{\Xcal}$}\end{matrix}\right\},\\
f &\lsends f^*
\end{aligned}
\end{equation}
is an equivalence.
\end{proposition}

\begin{proof}
As $S$ is quasi-compact, we find an fpqc-covering (or even a Zariski covering) of $S$ by an affine scheme. As both sides of \eqref{EqToniniTannaka} satisfy fpqc-descent in $S$, we may assume that $S = \Spec R$ is affine. 

The hypotheses on $S$ and $G$ imply that $\B{G}$ has an fpqc-covering by an affine scheme, hence $\B{G}$ is quasi-compact in the sense used in \cite{Tonini_Tannaka}. Moreover, $\B{G}$ has a quasi-affine diagonal because $G \to S$ is quasi-affine (Remark~\ref{DiagonalClassifying}). Hence we can apply \cite[Theorem~A]{Tonini_Tannaka} to see that \eqref{EqToniniTannaka} is an equivalence.
\end{proof} 

If $\Ycal$ is a noetherian algebraic stack, every quasi-coherent $\Oscr_{\Ycal}$-module is the colimit of its submodules of finite type by \cite[Proposition~15.4]{LMB_ChampsAlgebriques}). Hence $\Ycal$ has the resolution property if every coherent $\Oscr_{\Ycal}$-module is quotient of a vector bundle over $\Ycal$. Therefore \cite[\S 2]{Thomason_Equivariant} yields the following examples of classifying stacks with resolution properties.

\begin{example}\label{BGResolution}
Let $G$ be a group schemes over a scheme $S$. The classifying stack $\B{G}$ has the resolution property if one of the following conditions is satisfied.
\begin{assertionlist}
\item
The scheme $S$ is a separated regular noetherian scheme of dimension $\leq 1$ and $G \to S$ is an affine flat group scheme of finite type. 
\item
The scheme $S$ is a separated regular noetherian scheme of dimension $\leq 2$ and $G \to S$ is affine smooth with connected fibers.
\item
The scheme $S$ is a normal noetherian scheme and has an ample family of line bundles\footnote{Or, more generally, has itself the resolution property} and $G \to S$ is reductive.
\end{assertionlist}
\end{example}

For further examples, see, for instance, \cite{Gross_ResolutionPropertySurface}.


\subsection{$G$-bundles on classifying stacks}\label{Sec:GBundlesClassStack}

\begin{lemma}\label{MapsClassifyingStacks}
Let $H$ be a group object in a topos attached to some site $\Ccal$, and let $\Xcal$ be a stack on $\Ccal$. Then the groupoid of maps of stacks $\B{H} \to \Xcal$ has as objects pairs $(x \in \Hom(*,\Xcal), \varphi\colon H \to \Autline(x))$ and as morphisms the morphisms in $\Hom(*,\Xcal)$ that conjugate the maps of sheaves of groups $H \to \Autline(x)$.
\end{lemma}

\begin{proof}
A map $F\colon \B{H} \to \Xcal$ is sent to $(F(*), H = \Aut(*) \to \Aut(F(*)))$. This induces a bijection as is easily checked by identifying $\Hom(\B{H},\Xcal)$ with $\Hom([H\,{}_p\!\backslash *],\Xcal)$, where $[H\,{}_p\!\backslash *]$ is the prestack quotient.
\end{proof}

For $\Xcal = \B{G}$ this implies immediately the following more precise assertions.

\begin{remark}\label{HomClassifyingStack}
Let $G$ and $H$ be group algebraic spaces over $S$.
\begin{assertionlist}
\item
The groupoid $\Hom_S(\B{H},\B{G}) = \Bun_G(\B{H})$ is given by the groupoid consisting of a $G$-bundle $P$ over $S$ together with a homomorphism of group algebraic spaces $\mu\colon H \to \Autline_G(P)$. Here $(P,\mu)$ corresponds to the map
\[
\xymatrix{\B{H} \ar[r] & \B{\Autline_G(P)} \ar[rr]^-{\sim}_-{(\Theta^P)^{-1}} & & \B{G},}
\]
where the first map is induced by $\mu$ and the second is the inverse of the twist $\Theta^P$ by $P$ \cite[III, Section~2.3]{Giraud_NonAb}.

A morphism $(P,\mu) \to (P',\mu')$ is a morphism $u\colon P \to P'$ of $G$-bundles (necessarily an isomorphism) such that $\mu' = \textup{int}(u) \circ \mu$, where $\textup{int}(u)\colon \Autline(P) \to \Autline(P')$ is given by $\alpha \sends u \circ \alpha \circ u^{-1}$.
\item
Composition with the canonical map $S \to \B{H}$ defines an essentially surjective map of groupoids
\[
\gamma\colon \Hom_S(\B{H},\B{G}) \lto \Hom_S(S,\B{G}) = \Bun_G(S), \qquad (P,\alpha) \sends P.
\]
This construction also yields a map of stacks, again denoted by $\gamma$,
\[
\gamma\colon \Bunline_G(\B{H}) = \Homline_S(\B{H},\B{G}) \lto \Homline_S(S,\B{G}) = \B{G}.
\]
\item
Fix a $G$-bundle $P_0$ over $S$, corresponding to a map of stacks $S \to \B{G}$, and let $G_0 := \Autline_G(P_0)$ be the attached strong inner form of $G$. Then one has a 2-cartesian diagrams of stacks.
\begin{equation}\label{EqDescribeGBundleClassStack1}
\begin{aligned}\xymatrix{
\Homline_{\textup{$S$-Grp}}(H,G_0) \ar[r] \ar[d] & S \ar[d]^{P_0} \\
\Homline_S(\B{H},\B{G}) \ar[r]^-{\gamma} & \B{G}.
}\end{aligned}
\end{equation}
\item
Now suppose that $G \to S$ is an fpqc covering. Then $G_0$ is also an fpqc-covering of $S$. The morphism $P_0\colon S \to \B{G}$ can be viewed via twisting as a $G_0$-bundle. The left vertical arrow of \eqref{EqDescribeGBundleClassStack1} induces therefore an isomorphism
\begin{equation}\label{EqDescribeGBundleClassStack2}
[G_0\backslash \Homline_{\textup{$S$-Grp}}(H,G_0)] \liso \Bunline_G(\B{H}).
\end{equation}
\end{assertionlist}
\end{remark}

\noindent From the description \eqref{EqDescribeGBundleClassStack2} we obtain the following result using a result of Cotner.

\begin{proposition}\label{BunGClassStackAlgebraic}
Let $S$ be an algebraic space, let $G$ and $H$ be smooth affine group algebraic spaces over $S$. Suppose that the identity component $H^0$ of $H$ is a parabolic subgroup of a reductive group scheme over $S$ (e.g., if $H^0$ is reductive) and that the quotient $H^0\backslash H$ is finite \'etale over $S$.
\begin{assertionlist}
\item\label{BunGClassStackAlgebraic1}
The prestack $\Bunline_G(\B{H})$ is a stack for the fpqc topology and an algebraic stack locally of finite presentation over $S$.
\item\label{BunGClassStackAlgebraic2}
If $T \subseteq H^0$ is a maximal torus, then the restriction morphism $\Bunline_G(\B{H}) \to \Bunline_G(\B{T})$ is an affine morphism. If $H$ is reductive and $P \subseteq H$ is a parabolic subgroup, then $\Bunline_G(\B{H}) \to \Bunline_G(\B{P})$ is an open immersion.
\end{assertionlist}
\end{proposition}

\begin{proof}
By \cite[Theorem~5.8]{Cotner_HomSchemes}, $\Homline_{\textup{$S$-Grp}}(H,G)$ is an algebraic space which is ind-quasi-affine and locally of finite presentation over $S$. Hence \ref{BunGClassStackAlgebraic1} follows by \eqref{EqDescribeGBundleClassStack2}. 

Assertion~\ref{BunGClassStackAlgebraic2} then follows from \cite[Lemma~5.9, Lemma~5.10]{Cotner_HomSchemes}.
\end{proof}


\section{Restriction and extension of $G$-bundles}

%
%

In this section, $S$ denotes an algebraic space. All stacks are stacks for the fpqc-topology.


\subsection{Fully faithfulness of restriction}

\begin{proposition}\label{SectionAffine}
Let $\Xcal$ be an stack over $S$ whose diagonal is representable by algebraic spaces and that admits an fpqc covering by an algebraic space. Let $f\colon \Ucal \to \Xcal$ be a qcqs morphism of stacks that is representable by schemes. Suppose that $f_*\Oscr_{\Ucal} = \Oscr_{\Xcal}$.

Suppose we are given a 2-commutative diagram of solid arrows of stacks
\[\xymatrix{
\Ucal \ar[r] \ar[d]^f & \Ycal \ar[d]^g \\
\Xcal \ar[r] \ar@{.>}[ru] & \Zcal
}\]
with $g\colon \Ycal \to \Zcal$ a morphism representable by schemes. If $g$ is quasi-affine (resp.~affine), then there exists at most one (resp.~a unique)
morphism $\Xcal \to \Ycal$ making the above diagram commutative.
\end{proposition}

The condition that the diagonal of $\Xcal$ is representable by algebraic spaces ensures that every morphism $X \to \Xcal$, where $X$ is an algebraic space, is representable by algebraic spaces. In particular we know what it means for such a morphism to be an fpqc-covering.

\begin{proof}
Replacing $\Ycal$ by $\Ycal \times_{\Zcal} \Xcal$ and $g$ by its base change $\Ycal \times_{\Zcal} \Xcal \to \Xcal$, it suffices to show that for every representable quasi-affine (resp.~affine) morphism $g\colon \Ycal \to \Xcal$, every section of $g$ over $\Ucal$ has at most one (resp.~a unique) extension to $\Xcal$. We can work fpqc locally on $\Xcal$ and hence can assume that $\Xcal = \Spec A$ is an affine scheme. Then $f\colon \Ucal \to \Xcal$ is a qcqs map of schemes. Since $f_*\Oscr_{\Ucal} = \Oscr_{\Xcal}$, we have in particular $\Gamma(\Xcal,\Oscr_{\Xcal}) = \Gamma(\Ucal,\Oscr_{\Ucal})$.

Suppose that $g$ is affine. Then $\Ycal = \Spec B$ is an affine scheme and we find
\begin{align*}
\Hom_{\Xcal}(\Ucal,\Ycal) &= \Hom_A(B,\Gamma(\Ucal,\Oscr_{\Ucal})) = \Hom_A(B,A) \\
&= \Hom_{\Xcal}(\Xcal,\Ycal).
\end{align*}
Now suppose $g$ is quasi-affine and set $B := \Gamma(\Ycal,\Oscr_{\Ycal})$. The canonical map $\Ycal \to \Spec B$ is an open immersion and we have a commutative diagram
\[\xymatrix{
\Hom_{\Xcal}(\Xcal,\Ycal) \ar[r] \ar@{^{(}->}[d] & \Hom_{\Xcal}(\Ucal,\Ycal) \ar@{^{(}->}[d] \\
\Hom_{\Xcal}(\Xcal,\Spec B) \ar[r]^{\sim} & \Hom_{\Xcal}(\Ucal,\Spec B),
}\]
where the vertical maps are injective and the lower horizontal map is bijective. Therefore, the upper horizontal map is injective.
\end{proof}

Let us make some remarks about the hypotheses of the above proposition.

\begin{remark}\label{CheckFpqcLocally}
Let $\Xcal$ be a stack whose diagonal is representable by algebraic spaces and that admits an fpqc-covering by an algebraic space and let $f\colon \Ucal \to \Xcal$ be a qcqs morphism of stacks that is representable by schemes.
\begin{assertionlist}
\item
If $\Xcal$ is an affine scheme, then $\Oscr_{\Xcal} \to f_*\Oscr_{\Ucal}$ is an isomorphism if and only if $\Gamma(\Xcal,\Oscr_{\Xcal}) \to \Gamma(\Ucal,\Oscr_{\Ucal})$ is an isomorphism.
\item
To check that $\Oscr_{\Xcal} \to f_*\Oscr_{\Ucal}$ is an isomorphism, one can work fpqc-locally on $\Xcal$ since the formation of $f_*\Oscr_{\Ucal}$ commutes with flat base change.
\item
If $\Xcal$ is a locally noetherian scheme and $f\colon \Ucal \to \Xcal$ is an open immersion, then $\Oscr_{\Xcal} \to f_*\Oscr_{\Ucal}$ is an isomorphism if and only if the open subscheme $f(\Ucal)$ contains every point of depth $\leq 1$ of $\Xcal$ \cite[Th\'eor\`eme~5.10.5]{EGA4.2}.
\item
If $f$ is proper, flat, of finite presentation and has geometrically connected and geometrically reduced fibers, then $\Oscr_{\Xcal} = f_*\Oscr_{\Ucal}$.

Indeed, as we can work fpqc-locally on $\Xcal$, we can assume that $\Xcal$ is an algebraic space and even a scheme. Then we can conclude by \cite[Corollary~24.63]{GW2}.
\end{assertionlist} 
\end{remark}

\begin{corollary}\label{ExtendMapStackUniqueness}
Let $\Xcal$ be a stack over $S$ whose diagonal is representable by algebraic spaces and that admits an fpqc covering by an algebraic space, and let $f\colon \Ucal \to \Xcal$ be a morphism of stacks that is representable by schemes and qcqs. Suppose that $f_*\Oscr_{\Ucal} = \Oscr_{\Xcal}$. Let $\Ycal \to \Zcal$ be a morphism of stacks with representable affine diagonal.

Then composition with $f$ yields a fully faithful map of groupoids
\[
\Hom_{\Zcal}(\Xcal,\Ycal) \lto \Hom_{\Zcal}(\Ucal,\Ycal)
\]
\end{corollary}

\begin{proof}
Given any two $\Zcal$-morphisms $f,g \colon \Xcal \to \Ycal$, the morphism $\Isom_{\Zcal}(f,g) \to \Xcal$ is affine by hypothesis. Hence any section over $\Ucal$ extends uniquely to a section over $\Xcal$ by Proposition~\ref{SectionAffine}.
\end{proof}

\begin{proposition}\label{RestrictionFullyFaithful}
Let $G$ be a group algebraic space over $S$ such that $G \to S$ is flat and affine. Let $\Xcal$ be a stack over $S$ whose diagonal is representable by algebraic spaces and which admits an fpqc covering by an algebraic space, and let $f\colon \Ucal \to \Xcal$ be a qcqs morphism of stacks over $S$ that is representable by schemes. Suppose that $\Oscr_{\Xcal} \to f_*\Oscr_{\Ucal}$ is an isomorphism.

Then the functors
\begin{align}
f^*\colon \Vec{\Xcal} &\to \Vec{\Ucal},\label{EqRestrictionVec}\\
f^*\colon \Bun_G(\Xcal) &\to \Bun_G(\Ucal)\label{EqRestrictionGBun}
\end{align}
are fully faithful. 
\end{proposition}

\begin{proof}
It is standard that \eqref{EqRestrictionVec} is fully faithful which can be proved as for morphism of schemes (e.g. \cite[Lemma~24.67]{GW2}).

The fully faithfulness of \eqref{EqRestrictionGBun} follows from Corollary~\ref{ExtendMapStackUniqueness} applied to $\Ycal = \B{G}$ since $\B{G}$ has an affine diagonal by Remark~\ref{DiagonalClassifying}.
\end{proof}


\subsection{Extension of $G$-bundles}

Let us give criteria when a vector bundle (resp.~a $G$-bundle) is in the essential image of \eqref{EqRestrictionVec} (resp.~\eqref{EqRestrictionGBun}).

\begin{proposition}\label{InEssentialImage}
Suppose the hypotheses of Proposition~\ref{RestrictionFullyFaithful} are satisfied.
\begin{assertionlist}
\item\label{InEssentialImage1}
The question whether a vector bundle or a $G$-bundle on $\Ucal$ is in the essential image of $f^*$ can be answered fpqc-locally on $\Xcal$ since we already know that $f^*$ is fully faithful after any base change by an fpqc morphism.
\item\label{InEssentialImage2}
Under the hypotheses of the theorem, a vector bundle $\Fscr$ on $\Ucal$ is in the essential image of $f^*$ if and only if the direct image $f_*\Fscr$ is a vector bundle and the canonical map $f^*f_*\Fscr \to \Fscr$ is an isomorphism.
\item\label{InEssentialImage3}
A $G$-bundle $P$ on $\Ucal$ is in the essential image of \eqref{EqRestrictionGBun} if (and only if) locally for the fpqc-topology on $S$ there exists a monomorphism $\iota\colon G \mono H$ of relatively affine and flat group algebraic spaces over $S$ such that the fpqc-quotient $G\backslash H$ is representable by an algebraic space that is affine over $S$ such that $\iota_*(P)$ is in the essential image of $f^*\colon \Bun_H(\Xcal) \to \Bun_H(\Ucal)$.
\end{assertionlist}
\end{proposition}

\begin{proof}
Assertion~\ref{InEssentialImage1} follows formally from that fact that $f^*$ is fully faithful and that fpqc-descent is effective for vector bundles and for $G$-bundles. Assertion~\ref{InEssentialImage2} can be shown as in as in \cite[Lemma~24.67]{GW2} where the argument is given for a morphism of schemes.

Let us show \ref{InEssentialImage3}. By \ref{InEssentialImage1} we may assume that $\iota$ exists globally. Consider the map $\bar\iota\colon \B{G} \to \B{H}$ induced by $\iota$. Since $\iota$ is a monomorphism and $G\backslash H$ is affine, $\bar\iota$ is affine. The condition translates into the existence of a 2-commutative diagram
\[\xymatrix{
\Ucal \ar[r]^-P \ar[d]_f & \B{G} \ar[d]^{\bar\iota} \\
\Xcal \ar[r] & \B{H}.
}\]
Now Proposition~\ref{SectionAffine} yields a map $\Xcal \to \B{G}$ commuting with the diagram, i.e. a $G$-bundle on $\Xcal$ extending $P$.\end{proof}

In practice, we will usually apply Criterion~\ref{InEssentialImage3} of Proposition~\ref{InEssentialImage} with $H = \GL_m$ to reduce the essential surjectivity to a problem of vector bundles. Thus to apply this criterion we want to give examples for $G$ to satisfy the following properties.
\begin{assertionlist}
\item\label{EmbCrit1}
There exists an embedding $G \to \GL_m$ fpqc-locally.
\item\label{EmbCrit2}
Given an embedding $G \to \GL_m$, the quotient $G\backslash \GL_m$ is affine.
\end{assertionlist}

With respect to~\ref{EmbCrit1} there are the following results.

\begin{remark}\label{GeometricallyReductiveEmbeddable}
\begin{assertionlist}
\item\label{GeometricallyReductiveEmbeddable1}
Let $S$ be a noetherian scheme and let $G \to S$ be flat, affine and of finite type. Suppose that $S$ is regular of dimension $\leq 2$ \emph{or} that $\B{G}$ satisfies the resolution property. Then there exists a closed immersion $G \mono \GL(\Escr)$ for a vector bundle $\Escr$ on $S$ by \cite[Exp.~6B, Proposition~13.2 and Proposition~13.5]{SGA3NeuI}. 

If $S$ is affine, then $\Escr$ is a direct summand of a free vector bundle and one obtains a closed immersion $G \mono \GL_m$ for some $m$.
\item\label{GeometricallyReductiveEmbeddable2}
Let $G$ be a reductive group scheme over an algebraic space $S$. Then there exists \'etale locally on $S$ a closed immersion $G \mono \GL_m$ for some $m$.

Indeed, locally for the \'etale topology there exists a splitting of $G$ and a split reductive group is obtained by base change from a reductive group over $\ZZ$ which admits a closed immersion into some $\GL_m$ by \ref{GeometricallyReductiveEmbeddable1}.
\end{assertionlist}
\end{remark}

To deal with the question whether $G\backslash \GL_m$ is affine, recall the following notion, first studied by Franjou and van der Kallen in \cite{FranjouVanDerKallen} under the notion ``power-reductive''.

\begin{definition}\label{DefGeomReductive}
Let $S$ be an algebraic space. A flat, separated, finitely presented group algebraic space $G \to S$ is called \emph{geometrically reductive} if for every surjective map of quasi-coherent $G$-$\Oscr_S$-algebras $\Ascr \to \Bscr$ the map $\Ascr^G \to \Bscr^G$ is adequate, i.e., for every local section $s$ of $\Bscr^G$, \'etale locally some power of $s$ is in the image.
\end{definition}

Here we follow the terminology of Alper introduced in \cite{Alper_Adequate}. It strengthens the notion of geometric reductivity defined by Seshadri \cite{Seshadri_GeomRedArbBases}.

The property of geometric reductivity is stable under base change and under fpqc-descent \cite[Proposition~9.3.1]{Alper_Adequate}.

For us this notion is useful because of the following results.

\begin{remark}\label{RemGeomRedGrp}
Let $S$ be an algebraic space and let $G \to S$ be a flat, separated, finitely presented group algebraic space.
\begin{assertionlist}
\item\label{RemGeomRedGrp0}
If $G$ is of multiplicative type, then taking invariants is an exact functor, hence $G$ is geometrically reductive.
\item\label{RemGeomRedGrp1}
If $G$ a reductive group scheme over $S$, then $G$ is geometrically reductive \cite[Theorem~12]{FranjouVanDerKallen}. More generally, if $G$ is a smooth affine group scheme over $S$ such that its identity component $G^0$ is reductive and such that $G/G^0$ is finite, then $G$ is geometrically reductive \cite[Theorem~9.7.6]{Alper_Adequate}.
\item\label{RemGeomRedGrp2}
Let $H$ be an affine geometrically reductive group over $S$ and let $G \mono H$ be a monomorphism. Then the quotient $G\backslash H$ is affine if and only if $G$ is geometrically reductive \cite[Theorem~9.4.1]{Alper_Adequate}.
\end{assertionlist}
\end{remark}

Thus we can apply Proposition~\ref{InEssentialImage} and \ref{RestrictionFullyFaithful} as follows which is used in the classification of $G$-bundles on the projective line below (Theorem~\ref{ClassifyGBundlesPP1}).

\begin{proposition}\label{ExtendGBundleOpenImmersion}
Let $S$ be an algebraic space and let $D$ be a group algebraic space over $S$ such that $D \to S$ is ind-quasi-affine and an fpqc covering. Let $X$ be a scheme over $S$ endowed with a $D$-action and let $j\colon U \mono X$ be a $D$-invariant open subscheme such that $\Oscr_X = j_*\Oscr_U$. Let $G$ be a group algebraic space over $S$ such that $G \to S$ is flat and affine.
\begin{assertionlist}
\item\label{ExtendGBundleOpenImmersion1}
The restriction yields a fully faithful functor
\begin{equation}\label{EqResQuotientStack}
\Bun_G([D\backslash X]) \lto \Bun_G([D\backslash U]).
\end{equation}
\item\label{ExtendGBundleOpenImmersion2}
Suppose that $S$ and $X$ are noetherian regular of dimension $\leq 2$ and that $X \setminus U$ is of codimension $2$ in $X$. If $G$ is an extension of finite \'etale group scheme by a reductive group scheme, then \eqref{EqResQuotientStack} is an equivalence.
\end{assertionlist}
\end{proposition}

\begin{proof}
Set $\Xcal := [D\backslash X]$ and $\Ucal := [D\backslash U]$ and let $j\colon \Ucal \to \Xcal$ be the inclusion. Let us show~\ref{ExtendGBundleOpenImmersion1}. We have $\Oscr_{\Xcal} = j_*\Oscr_{\Ucal}$ by Remark~\ref{CheckFpqcLocally} and thus conclude by Proposition~\ref{RestrictionFullyFaithful}.

Let us show \ref{ExtendGBundleOpenImmersion2}. Since we can show essential surjectivity fpqc-locally on $[D \backslash X]$ by Remark~\ref{InEssentialImage}~\ref{InEssentialImage1}, it suffices to show that the fully faithful functor $\Bun_G(X) \to \Bun_G(U)$ is essentially surjective. Moreover, we can assume that $S$ and $X$ are affine. Then $G$ admits an embedding into some $\GL_m$ by Remark~\ref{GeometricallyReductiveEmbeddable}~\ref{GeometricallyReductiveEmbeddable1}. As $G$ is geometrically reductive by Remark~\ref{RemGeomRedGrp}~\ref{RemGeomRedGrp1}, we see that $G\backslash\GL_m$ is affine by Remark~\ref{RemGeomRedGrp}~\ref{RemGeomRedGrp2}. Thus we may assume $G = \GL_m$ by Proposition~\ref{InEssentialImage}~\ref{InEssentialImage3} and we are left to prove that every vector bundle on $U$ extends to $X$. But this is a classical result from commutative algebra (e.g.~\cite[0EBJ, 0B3N]{Stacks}).
\end{proof}

Let us give one further application of the results above.

\begin{proposition}\label{PullBackProper}
Let $m \geq 1 $ be an integer. Let $S$ be an algebraic space and let $G$ be a geometrically reductive group algebraic space over $S$ that admits locally for the \'etale topology a closed embedding into $\GL(\Escr)$ for some vector bundle $\Escr$ on $S$ (for instance, if $G$ is a reductive group scheme over $S$, see Remark~\ref{GeometricallyReductiveEmbeddable}~\ref{GeometricallyReductiveEmbeddable2}).

Let $\Xcal$ be a stack over $S$ whose diagonal is representable by algebraic spaces and that admits an fpqc covering by an algebraic space. Let $f\colon \Ucal \to \Xcal$ be a morphism of stacks that is representable by schemes and that is proper, flat, of finite presentation with geometrically connected and geometrically reduced fibers. Then
\[
f^*\colon \Bun_G(\Xcal) \to \Bun_G(\Ucal) \tag{*}
\]
is fully faithful.

Let $P$ be a $G$-bundle on $\Ucal$ and suppose there exists a reduced scheme $Y$ and an fpqc-covering $Y \to \Xcal$ such that the pullback of $P$ to $V := \Ucal \times_{\Xcal} Y$ is trivial on the fibers of $V \to Y$. Then $P$ is in the essential image of (*).
\end{proposition}

\begin{proof}
The hypothesis on $G$ implies that $G \to S$ is affine. Let us show that $\Oscr_{\Xcal} = f_*\Oscr_{\Ucal}$. We can do this fpqc-locally and therefore assume that $\Xcal$ is an algebraic space or even a scheme. In this case, it is well known that the hypotheses on $f$ imply $\Oscr_{\Xcal} = f_*\Oscr_{\Ucal}$ (e.g., \cite[Corollary~24.63]{GW2}). Therefore, the functor $f^*$ is fully faithful by Proposition~\ref{RestrictionFullyFaithful}.

Let $S' \to S$ be an \'etale covering such that $G \times_S S'$ admits a closed embedding into $\GL_{m,S'}$. Let $\Xcal' := \Xcal \times_S S'$ be the base change of $\Xcal$ and define similarly $\Ucal'$, $f'$, $Y'$, and $V'$. By Proposition~\ref{InEssentialImage}~\ref{InEssentialImage1}, it suffices to show that the pullback of $P$ to $\Ucal'$ is in the image of $f^{\prime*}\colon \Bun_G(\Xcal') \to \Bun_G(\Ucal')$. Since $Y' \to Y$ is \'etale, $Y'$ is also reduced. Moreover, the pullback of $P$ to $V'$ is also trivial on all fibers of $V' \to Y'$. Therefore, we can assume that $G$ admits a closed embedding into $\GL_{m,S}$. As $G$ is geometrically reductive, $G\backslash \GL_m$ is affine.

By Proposition~\ref{InEssentialImage}~\ref{InEssentialImage3} it therefore suffices to show that every vector bundle $\Escr$ on $V$ which is trivial on all fibers of $V \to Y$ is isomorphic to the pullback of a vector bundle on $Y$. Now the Seesaw Theorem in the form of \cite[Theorem~24.66]{GW2} yields a locally closed subscheme $Z$ of $Y$ whose underlying topological space is equal to $Y$ and a vector bundle on $Z$ whose pullback is the restriction of $\Escr$ to $V \times_Y Z$. Since $Y$ is reduced, we must have $Z = Y$.
\end{proof}

%% file: Ch2GBundlesProjectiveLine.tex
\section{Lifting of $G$-bundles}

%
%
%

In this section, the main result is a lifting result for $G$-bundles for certain henselian pairs of algebraic stacks (Theorem~\ref{LiftHenselian}). We start with a purely topological result that allows to describe the henselian property via spaces of connected components (Proposition~\ref{PropHenselianMorph}).

\subsection{Connected components and clopen subsets}\label{Sec:ClopenConnComp}

For a topological space $X$ we denote by $\pi_0(X)$ its space of connected components, i.e., the underlying set is the set of connected components of $X$ and $\pi_0(X)$ is endowed with the quotient topology of $X$. This is a totally disconnected topological space. If $f\colon X \to Y$ is a continuous map of topological spaces, sending a connected component $C \subseteq X$ to the unique connected component of $Y$ that contains $f(C)$ yields a continuous map $\pi_0(f)\colon \pi_0(X) \to \pi_0(Y)$. We obtain a functor $\pi_0$ from the category of topological spaces to the category of totally disconnected spaces. This functor is left adjoint to the inclusion functor.


Let $X$ be a spectral topological space. We call a subset $Z$ of $X$ \emph{pro-clopen} if it is the intersection of subsets of $X$ that are both open and closed. By \cite[04PL, 0906]{Stacks}, the following assertions hold.

\begin{lemma}\label{LemmaT}
\begin{assertionlist}
\item\label{LemmaT2}
A subset is pro-clopen if and only if it is closed in $X$ and a union of connected components of $X$. 
\item\label{LemmaT3}
The space $\pi_0(X)$ is pro-finite and in particular compact Hausdorff.
\end{assertionlist}
\end{lemma}

From these results one easily deduces the following lemma.

\begin{lemma}\label{ConnPreSpectral}
Let $X$ be a spectral space. Let $p\colon X \to \pi_0(X)$ be the canonical quotient map. Then sending $Y \subseteq \pi_0(X)$ to $p^{-1}(Y) \subseteq X$ yields inclusion preserving bijections
\begin{align*}
\{\text{$Y \subseteq \pi_0(X)$ closed}\} &\iso \{\text{$Z \subseteq X$ pro-clopen}\}, \\
\{\text{$Y = \{C\} \subseteq \pi_0(X)$ singleton}\} &\iso \{\text{$\emptyset \ne Z \subseteq X$ minimal pro-clopen}\}, \\
\{\text{$Y \subseteq \pi_0(X)$ open and closed}\} &\iso \{\text{$Z \subseteq X$ open and closed}\},
\end{align*}
The inverse map is given by $Z \sends p(Z)$.
\end{lemma}

It follows that the minimal non-empty pro-clopen subsets of $X$ are precisely the connected components of $X$.

\begin{proof}
A subset $Z$ of $X$ is of the form $p^{-1}(Y)$ for some subset $Y$ of $\pi_0(X)$ if and only if $Z$ is a union of connected components of $X$. In this case $Z = p^{-1}(Y)$ is open (resp.~closed) in $X$ if and only if $Y$ is open (resp.~closed) in $\pi_0(X)$. Hence the first bijection follows from Lemma~\ref{LemmaT}~\ref{LemmaT2}. The second bijection follows because in $\pi_0(X)$ the minimal non-empty closed subsets are the singletons since $\pi_0(X)$ is Hausdorff. The third bijection follows since any open and closed subset $Z$ of $X$ is a union of connected components and hence of the form $p^{-1}(Y)$ for some open and closed subset $Y$ of $\pi_0(X)$.  
\end{proof}

For a topological space $X$ we denote by $\ClOpen(X)$ the set of open and closed subsets of $X$. 

\begin{proposition}\label{ConnMapPreSpektral}
Let $f\colon X \to Y$ be a continuous map of spectral topological spaces. Then the following assertions are equivalent.
\begin{equivlist}
\item\label{ConnMapPreSpektral1}
Sending $V \subseteq Y$ to $f^{-1}(V) \subseteq X$ induces a bijection $\ClOpen(Y) \iso \ClOpen(X)$.
\item\label{ConnMapPreSpektral2}
The map $\pi_0(f)\colon \pi_0(X) \to \pi_0(Y)$ is bijective.
\item\label{ConnMapPreSpektral3}
The map $\pi_0(f)\colon \pi_0(X) \to \pi_0(Y)$ is a homeomorphism.
\end{equivlist}
\end{proposition}

These conditions are, for instance, satisfied if $f$ has connected fibers (so that in particular, $f$ is surjective) and if $f$ endows $Y$ with the quotient topology of $X$ because then $\pi_0(Y)$ satisfies the universal property of $\pi_0(X)$.

\begin{proof}
Since $\pi_0(X)$ and $\pi_0(Y)$ are compact Hausdorff, \ref{ConnMapPreSpektral2} and \ref{ConnMapPreSpektral3} are equivalent.

By the third bijection of Lemma~\ref{ConnPreSpectral}, Assertion~\ref{ConnMapPreSpektral1} is equivalent to the fact that $\pi_0(f)^{-1}(-)$ induces a bijection between the open and closed subsets of $\pi_0(Y)$ and of $\pi_0(X)$. This is clearly the case if $\pi_0(f)$ is a homeomorphism. This shows that \ref{ConnMapPreSpektral3} implies \ref{ConnMapPreSpektral1} and that the converse follows from the following lemma.
%
%
\end{proof}

\begin{lemma}\label{HomeoProfinite}
Let $f\colon X \to Y$ be a continuous map of pro-finite spaces such that $f^{-1}(-)$ induces a bijection $\ClOpen(Y) \iso \ClOpen(X)$. Then $f$ is a homeomorphism.
\end{lemma}

\begin{proof}
It suffices to show that $f$ is bijective. Let us show that $f$ is injective. Let $x_1,x_2 \in X$ with $f(x_1) = f(x_2) =: y$. Let $A$ be an open and closed neighborhood of $x_1$. By hypothesis, there exists an open and closed neighborhood $B \subseteq Y$ of $y$ with $f^{-1}(B) = A$. Hence $A$ also contains $x_2$. As each singleton in a pro-finite space is the intersection of its open and closed neighborhoods, we see $x_1 = x_2$.

It remains to show that $f$ is surjective. Let $y \in Y$ and write $\{y\} = \bigcap_{y\in B}B$ as the intersection of the open and closed neighborhoods $B$ of $y$. The $B$'s form a cofiltered system and $f^{-1}(B)$ is non-empty by hypothesis since $B$ is non-empty. Hence $f^{-1}(y) = \bigcap_{y\in B}f^{-1}(B) \ne \emptyset$ because $X$ is quasi-compact.
\end{proof}


\subsection{Henselian pairs of stacks}

We call an fpqc-stack $\Xcal$ \emph{algebraic} if its diagonal is representable by algebraic spaces and if there exists a surjective and smooth morphism $U \to \Xcal$, where $U$ is a scheme.

In the stacks project it is only assumed that $\Xcal$ is an fppf-stack. Since our main results use only algebraic stacks that are quotient stacks by affine flat finitely presented group schemes $G$ and since all $G$-bundles for the fpqc-topology are automatically fppf-locally trivial, this difference in terminology will not trouble us.

%
We denote by $|\Xcal|$ the underlying topological space of an algebraic stack $\Xcal$ \cite[04Y8]{Stacks}.

For any algebraic stack $\Xcal$, we denote by $\ClOpen(\Xcal)$ the set of open and closed substacks of $\Xcal$. We denote by $\pi_0(\Xcal)$ the space of connected components of the underlying topological space of $\Xcal$. If $\Xcal$ is qcqs, the underlying topological space of $\Xcal$ is spectral \cite[0DQP]{Stacks} and hence $\pi_0(\Xcal)$ is pro-finite \cite[0906]{Stacks}.

Recall the notion of a henselian pair for rings.

\begin{defrem}\label{DefHenselianPairRing}
Now let $A$ be a ring and $Z \subseteq X := \Spec A$ be a closed subscheme. Recall that $(X,Z)$ is called henselian if for every finite morphism $X' \to X$, the map $\ClOpen(X') \to \ClOpen(X' \times_X Z)$ is bijective. Clearly, the property for $(X,Z)$ to be henselian depends only on the underlying topological subspace of $Z$. There is a unique smallest closed subspace $Z_{\min}$ in $X$ such that $(X,Z_{\min})$ is henselian \cite[0G1S]{Stacks}. If $I \subseteq A$ is an ideal, one also says that $(A,I)$ is henselian if $(X,V(I))$ is henselian. See \cite[Theorem~20.15]{GW2} for other characterizations of henselian pairs $(A,I)$. 

Examples of henselian pairs are pairs $(A,I)$, where $I$ is an ideal such that $A$ is $I$-adically complete (or, more general, if $A$ is derived $I$-adically complete \cite[0G3H]{Stacks}).
\end{defrem}

In \cite[Theorem~2.1.6]{BC_Torsors} it is shown that for an affine henselian pair $(X,Z)$ as above and a quasi-affine smooth group scheme $G$ over $X$ restriction of $G$-bundles defines a bijective map
\[
H^1(X,G) \liso H^1(Z,G).
\]
We want to generalize this result to certain pairs $(\Xcal,\Zcal)$ consisting of an algebraic stack $\Xcal$ and a closed substack $\Zcal$.

We have the following result which follows from the purely topological statement Proposition~\ref{ConnMapPreSpektral}.

\begin{proposition}\label{PropHenselianMorph}
Let $h\colon \Zcal \to \Xcal$ be a representable morphism of qcqs algebraic stacks. For a morphism $\Xcal' \to \Xcal$ of algebraic stacks denote by $h'\colon \Zcal' := \Zcal \times_{\Xcal} \Xcal' \to \Xcal'$ the base change of $h$. Then the following conditions are equivalent.
\begin{equivlist}
\item\label{PropHenselianMorph1}
For every finite morphism $\Xcal' \to \Xcal$, the map
\[
\ClOpen(\Xcal') \lto \ClOpen(\Zcal'), \qquad U \sends h^{\prime-1}(U)
\]
is bijective.
\item\label{PropHenselianMorph2}
For every finite morphism $\Xcal' \to \Xcal$, the map
\[
\pi_0(\Zcal') \lto \pi_0(\Xcal')
\]
is a homeomorphism.
\end{equivlist}
\end{proposition}



As a special case one has the notion of a henselian pair of algebraic stacks:

\begin{definition}\label{DefHenselianPair}
A pair consisting of qcqs algebraic stack $\Xcal$ and a closed substack $\Zcal$ is called \emph{henselian}, if the inclusion $\Zcal \to \Xcal$ satisfies the equivalent conditions of Proposition~\ref{PropHenselianMorph}.
\end{definition}

%

The only source of henselian pairs of algebraic stacks that we will need in this paper are algebraic stacks that admit moduli spaces which form an henselian pair as in the following remark.

\begin{remark}\label{RemarkAdequateModuli}
Let $S$ be an algebraic space. Let $G$ be a geometrically reductive group scheme over $S$ (see Definition~\ref{DefGeomReductive}). Let $p\colon \Xtilde \to S$ be an affine morphism of algebraic spaces and let $\Xtilde$ be endowed with an action by $G$ over $S$. By \cite[Theorem~9.1.4]{Alper_Adequate} and \cite[Theorem~3.12]{AHR_EtaleLocal} the map
\[
\pi\colon [G\backslash \Xtilde] \lto \Spec(p_*\Oscr_{\Xtilde})^G
\]
is initial in the category of all maps from $[G\backslash \Xtilde]$ to algebraic spaces. The morphism $\pi$ is surjective and universally closed by \cite[Theorem~5.3.1]{Alper_Adequate}. The same reference also shows that the underlying topological space of each geometric fiber of $\pi$ contains at most one closed point. Moreover, since $\pi$ is qcqs, every geometric fiber is spectral. In particular, it is a quasi-compact Kolmogorov space and therefore always has at least one closed point \cite[005E]{Stacks}. Therefore every geometric fiber is spectral with a unique closed point and hence homeomorphic to the spectrum of a local ring. In particular, $\pi$ has geometrically connected fibers.

Let $\Zcal$ be a closed substack of $\Xcal := [G\backslash X]$ and let $Z$ be a closed subscheme of $X := \Spec(p_*\Oscr_{\Xtilde})^G$ such that $\pi(\Zcal) = Z$ (set-theoretically). Then $(\Xcal,\Zcal)$ is henselian if and only if $(X,Z)$ is henselian by \cite[Theorem~3.10]{AHR_EtaleLocal}.
\end{remark}

Below, we will also use the following remark.

\begin{remark}\label{StabilizerQuotientStack}
Let $G$ be a group algebraic space over $S$ such that $G \to S$ is an fpqc covering. Suppose $G$ acts on an algebraic space $X$ over $S$. Let $x\colon T \to [G\backslash X]$ be a $T$-valued point of the fpqc quotient where $T$ is an algebraic space, and let $G_x := \Autline_{\Xcal(T)}(x)$ be the fpqc-sheaf of automorphisms of $x$.

By definition of $[G\backslash X]$ as fpqc-stackification of $[G\,{}_p\!\backslash X]$, there exists an fpqc-covering $\pi\colon T' \to T$ and a morphism $\xtilde \colon T' \to X$ such that
\[\xymatrix{
T' \ar[r]^{\xtilde} \ar[d]_{\pi} & X \ar[d]^p \\
T \ar[r] & [G\backslash X]
}\]
commutes. Hence
\[
G_x \times_T T' = G_{x \circ \pi} = G_{p \circ \xtilde} = \Stab_{G_{T'}}(\xtilde),
\]
where $\Stab_{G_{T'}}(\xtilde)$ denotes the stabilizer of $\xtilde \in X(T')$ in $G_{T'}$. In particular, $G_x$ is always fpqc-locally on $T$ isomorphic to a subgroup algebraic space of $G_T$.
\end{remark}

The main result of this section is an application of \cite[Proposition~7.10]{AHR_EtaleLocal} that we use in the following situation.

\begin{void}\label{SetUpLifting}
Let $A$ be a ring and let $H$ be a group scheme over $A$ such that there exists a finite locally free surjective map $\Spec A' \to \Spec A$ such that $H \otimes_A A'$ is a diagonalizable group scheme of finite type. Let $B$ be an $A$-algebra and suppose that $H$ acts on $\Xtilde := \Spec B$. Let $\Ztilde = \Spec B/I$ be a closed $H$-invariant subscheme of $X$. Set $\Xcal := [H\backslash \Xtilde]$ and $\Zcal := [H\backslash \Ztilde]$ and let $X := \Spec B^H$ and $Z := \Spec (B/I)^H$. As $H$ is of multiplicative type, $Z$ is a closed subscheme of $X$ because $(-)^H$ is an exact functor.

We also fix a gerbe $\Gcal$ which is fpqc-locally isomorphic to $\B{G}$ for a flat, quasi-affine, finitely presented group scheme $G$ over $A$.
\end{void}

The condition on $H$ is, for instance, satisfied if $H$ is an arbitrary group scheme of multiplicative type of finite type over $A$ and $A$ is a normal domain: In that case $H$ even becomes diagonalizable after a finite \'etale surjective base change \cite[Corollary~B.3.6]{Conrad_Reductive}.

\begin{theorem}\label{LiftHenselian}
In the Situation~\ref{SetUpLifting} suppose that $(X,Z)$ is a henselian pair. Then any morphism $\Zcal \to \Gcal$ can be extended to a morphism $\Xcal \to \Gcal$. If $G$ is smooth\footnote{In \cite[Proposition~7.10]{AHR_EtaleLocal} the maybe misleading hypothesis is made that $\Gcal$ is a ``smooth gerbe'' over $A$ which might be interpreted as meaning that $\Gcal \to \Spec A$ is smooth. But this would include $\Gcal = \B{G}$ for a (not necessarily smooth) flat affine group scheme $G$ of finite presentation over $A$. But the proof of the uniqueness statement in \cite[Proposition~7.10]{AHR_EtaleLocal} uses that the cotangent complex of $\Gcal$ is perfect of tor-amplitude 1 (homological numbering), which does not hold for $\B{G}$ with $G$ non-smooth in general, see also Remark~\ref{SmoothnessNecessary} below.} and quasi-affine, this extension is unique up to non-unique 2-isomorphisms.
\end{theorem}

In particular, we see that for every smooth quasi-affine group scheme $G$ over $A$ the restriction functor
\[
\Bun_G(\Xcal) \lto \Bun_G(\Zcal)
\]
is essentially surjective and full. Therefore one has a bijection
\[
H^1(\Xcal,G) \liso H^1(\Zcal,G).\tag{*}
\]
If $G$ is flat, quasi-affine, and finitely presented over $A$, then (*) is still surjective.

\begin{proof}
We want to apply \cite[Proposition 7.10]{AHR_EtaleLocal} to the projection $\Gcal \times_{\Spec A} X \lto X$. For this we have to make sure that the hypotheses of that theorem are satisfied. As $H$ is of multiplicative type, $X$ and $Z$ are affine good moduli spaces of $\Xcal$ and $\Zcal$, respectively, by \cite[Theorem~13.2]{Alper_Good}. Moreover $\Xcal$ is clearly quasi-compact, and it has an affine diagonal by Proposition~\ref{QuotientByIndQuasiAffine}. It remains to check that the following hypotheses are satisfied.
\begin{definitionlist}
\item\label{LiftHenseliana}
The projection $\Gcal \times_{\Spec A} X \lto X$ is quasi-separated and smooth with affine stabilizers.
\item\label{LiftHenselianb}
Every quasi-coherent $\Oscr_{\Xcal}$-module of finite type is a quotient of a vector bundle on $\Xcal$,
\item\label{LiftHenselianc}
The stabilizer of every (closed) point of $\Xcal$ is of multiplicative type.
\end{definitionlist}
Let us show \ref{LiftHenseliana}. It suffices to show that $f\colon \B{G} \to \Spec A$ is quasi-separated and smooth with affine stabilizers. Since $G$ is quasi-compact and quasi-separated, $f$ has a quasi-compact and quasi-separated diagonal which means that $f$ is quasi-separated. The composition $\Spec A \to \B{G} \to \Spec A$ is the identity and in particular smooth. Since $\Spec A \to \B{G}$ is faithfully flat and of finite presentation, $f$ is smooth. Since $G$ is quasi-affine, $f$ has quasi-affine stabilizers. But any quasi-affine group scheme of finite type over a field is affine by \cite[$\textup{VI}_B$, Proposition~11.11]{SGA3I}. 

Every diagonalizable group scheme of finite type can be embedded into $\GG_m^n$ and in particular into $\GL_n$ for some $n \geq 1$. Therefore $H$ can be embedded into $\GL(E)$, where $E$ is the finite projective $A$-module $(A')^n$. Therefore Condition~\ref{LiftHenselianb} is satisfied by \cite[Remark~2.5]{AHR_EtaleLocal}.

Finally, for a field $k$ and a $k$-valued point $z\colon \Spec k \to \Zcal$ let $H_z$ be its stabilizer. By Example~\ref{StabilizerQuotientStack} there exists a finite extension $k \to k'$ such that $H_z \otimes_k k'$ is isomorphic to a subgroup scheme of $H \otimes_A k'$. Since every subgroup scheme of a group scheme of multiplicative type over a field is again of multiplicative type \cite[Exp.~IX, Proposition~8.1]{SGA3II}, the stabilizer $H_z$ is of multiplicative type.
\end{proof}

\begin{remark}\label{HypthosesisN}
The hypothesis that $H$ is diagonalizable after a finite locally free base change is not the only possible hypothesis for Theorem~\ref{LiftHenselian}. By using one of the other hypotheses of \cite[Proposition~7.10]{AHR_EtaleLocal} it would also suffice to assume that the following two conditions are satisfied instead of the condition on $H$.
\begin{definitionlist}
\item
Every point of $Z$ has positive characteristic {\em or} there are only a finite number of different characteristics of points of $Z$.
\item
There exists a closed embedding $H \to \GL_{n,A}$ for some $n \geq 1$ (this is, for instance, the case if $A$ is normal and noetherian by \cite[Corollary~9.10]{AHR_EtaleLocal}) and the group scheme $H$ is linearly reductive in the sense of \cite[Definition 12.1]{Alper_Good}.
\end{definitionlist}
\end{remark}


\subsection{Reformulation of the lifting theorem in terms of graded algebras}\label{SecHenselianGraded}

Let $M$ be an abelian group and let $H := \Homline(\Mline,\GG_{m})$ be the corresponding diagonalizable group scheme over $\ZZ$. Recall that $H$ is a flat and affine group scheme over $\ZZ$ \cite[VIII, Proposition~2.1]{SGA3II}. In this case, we can make Situation~\ref{SetUpLifting} more concrete as follows.

Let $A$ be a ring. By \cite[Exp.~I, Corollaire~4.7.3.1]{SGA3I} the functor $B \sends \Spec B$ induces a contravariant equivalence of the category of $M$-graded $A$-algebras and the category of affine $A$-schemes with a $H_A$-action. Hence in the sequel we will identify affine $A$-schemes with a $H_A$-action and $M$-graded $A$-algebras.

For an $M$-graded $A$-algebra $B$ we form the quotient stack $\Xcal := [H_A\backslash \Spec B]$. Let $I = \bigoplus_{m}I_m \subseteq B$ be a graded ideal and let $\Zcal := [H_A\backslash(\Spec B/I)]$ be the corresponding closed substack of $\Xcal$.


We obtain the following special case of Theorem~\ref{LiftHenselian}.

\begin{corollary}\label{LiftHenselianGraded}
With the above notation of suppose that $M$ is finitely generated and that $(B_0,I_0)$ is a henselian pair. Let $G$ be a flat, quasi-affine, finitely presented (resp.~smooth quasi-affine) group scheme over $A$.

Then the restriction functor $\Bun_G(\Xcal) \to \Bun_G(\Zcal)$ is essentially surjective (resp.~essentially surjective and full). In particular,
\[
H^1(\Xcal,G) \lto H^1(\Zcal,G)
\]
is surjective (resp.~bijective).
\end{corollary}

\begin{example}\label{TrivialHenselian}
As a special case suppose that $M = \ZZ$ and that $B$ is concentrated in degrees $\geq 0$. We choose $I := \bigoplus_{m>0}B_m$. Then $I_0 = 0$ and $(B_0,0)$ is trivially henselian. Hence one has for any smooth quasi-affine group scheme over $A$ a bijection
\begin{equation}\label{EqH1TrivialHensel}
H^1([\GG_m\backslash (\Spec B)],G) \iso H^1(\B{\GG_{m,B_0}},G).
\end{equation}
In this case, the maps $B_0 \to B \to B/I = B_0$ induce maps
\[
H^1(\B{\GG_{m,B_0}},G) \to H^1([(\Spec B)/\GG_m],G) \iso H^1(\B{\GG_{m,B_0}},G)
\]
whose composition is the identity. Here the first map is given by pullback along $\pi\colon [(\Spec B)/\GG_m] \to \B{\GG_{m,B_0}}$ induced by the structure map of $\Spec B$ as a $B_0$-scheme. This shows that the inverse of \eqref{EqH1TrivialHensel} is given by $\pi^*$.

The same assertions also holds if $B$ is concentrated in degrees $\leq 0$.
\end{example}


The above example is in fact a special case of considering the fixed point locus:

\begin{remark}\label{DescribeFixGGm}
Suppose that $A$, $B$, $X$ and $\Xcal$ are defined as above.

We claim that the fixed point locus of the $H_A$-scheme $X$ is given by the closed subscheme $\Spec (B/I^0)$, where $I^0$ is the graded ideal of $B$ generated by the homogenous elements in non-zero degree. Indeed, if $R$ is any $A$-algebra, then $X^{H_A}(R)$ is the subset of $M$-graded $A$-algebra homomorphisms $B \to R$ with $R$ trivially $M$-graded which shows the claim.

The inclusion $B_0 \to B$ induces an isomorphism of trivially $M$-graded $A$-algebras
\[
B^0 := B_0/I_0^0 \iso B/I^0,
\]
where $I^0_0 = \sum_{m\ne 0}B_{-m}B_m$ is the degree $0$ part of $I^0$.

The closed immersion $\Spec B^0 = (\Spec B)^{H_A} \mono \Spec B$ induces a closed immersion
\[
i\colon \B{H_{B^0}} = [H_A\backslash \Spec B^0] \mono [H_A\backslash \Spec B].
\]
If $H$ is of finite type over $\ZZ$, the pair $([H_A\backslash \Spec B], \B{H_{B^0}})$ is henselian if and only if $(B_0,I^0_0)$ is henselian. If this is the case, then $i^*$ yields a bijection
\begin{equation}\label{EqGBundleFixedPoint}
H^1(\Xcal,G) \iso H^1(\B{H_{B^0}},G)
\end{equation}
for every smooth quasi-affine group scheme $G$ over $A$. We refer to Remark~\ref{HomClassifyingStack} for a description of the right hand side of \eqref{EqGBundleFixedPoint}.
\end{remark}


\subsection{Lifting of vector bundles over quotient stacks by diagonalizable groups}

%

In this section we give an elementary proof of (a generalization of) Theorem~\ref{LiftHenselian} in the case of vector bundles and that $H$ is diagonalizable (but not necessary of finite type) with character group $M$.

Let us describe quasi-coherent modules and vector bundles on quotient stacks by $H$.

\begin{remdef}\label{ModulesQuotientGGm}
Let $A$ be a ring, let $B$ be an $M$-graded $A$-algebra and set $\Xcal := [H_A\backslash \Spec B]$. By descent along $\Spec B \to \Xcal$, quasi-coherent $\Oscr_{\Xcal}$-modules are identified with $H_A$-equivariant quasi-coherent modules on $\Spec B$ and these can be identified with $M$-graded $B$-modules. Morphisms are $B$-linear maps preserving the grading.

The tensor product of two $M$-graded $B$-modules $E$ and $F$ is the graded $B$-module whose underlying module is the tensor product $E \otimes_B F$ obtained by forgetting the gradings. It is endowed with the $M$-grading such that $(E \otimes_B F)_m$ is generated by $e \otimes f$ with $e \in E_k$, $f \in F_{l}$ such that $k + l = m$.

We obtain an monoidal equivalence of the symmetric monoidal abelian category of quasi-coherent modules on $\Xcal$ and the symmetric monoidal abelian category of $M$-graded $B$-modules. Via this equivalence, quasi-coherent $\Oscr_{\Xcal}$-modules of finite type correspond to $M$-graded $B$-modules that are finitely generated as $M$-graded $B$-modules (or, equivalently, that are finitely generated as $B$-modules).

Recall that an $M$-graded $B$-module $E$ is called \emph{projective} if it satisfies the following equivalent conditions.
\begin{equivlist}
\item
$E$ is projective as an object in the abelian category of $M$-graded $B$-modules.
\item
$E$ is projective as a $B$-module.
\item
$E$ is a direct summand of an $B$-module of the form $\bigoplus_i B(d_i)$, where $B(d_i)$ is the module $B$ endowed with the grading shifted by $d_i \in M$, i.e., $B(d_i)_n = B_{d_i+n}$.
\end{equivlist}
The category $\Vec{\Xcal}$ of vector bundles on $\Xcal$ is identified with the exact tensor category of finitely generated, projective, $M$-graded $B$-modules. 
\end{remdef}

\begin{example}\label{VectorbundlesoverAA1ModGGm}
Let $A$ be a ring. Then the category of quasi-coherent modules over $[\GG_{m,A}\backslash \AA_A^1]$ is equivalent to the category $\ZZ$-graded $A$-modules $E = \bigoplus_i E_i$ together with a graded $A$-linear endomorphism $T\colon E \to E$ of degree $1$ given by the indeterminate $T \in A[T] = \Oscr_{\AA^1_A}(\AA^1_A)$.

Such a quasi-coherent module $(E,T)$ corresponds to a vector bundle over $[\GG_{m,A}\backslash \AA_A^1]$ if and only if
\begin{definitionlist}
\item
$E_i$ is a finite projective $A$-module for all $i$,
\item
$T\colon E_i \to E_{i+1}$ is injective with projective cokernel,
\item
$E_i = 0$ for $i \ll 0$ and $T\colon E_i \to E_{i+1}$ is an isomorphism for $i \gg 0$.
\end{definitionlist}
In other word, setting
\[
E_{\infty} := \colim (\cdots \ltoover{T} E_i \ltoover{T} E_{i+1} \ltoover{T} \cdots ),
\]
a vector bundle over $[\GG_{m,A}\backslash \AA_A^1]$ is the same as a finite projective $A$-module $E_{\infty}$ together with a filtration by direct summands.
\end{example}

We continue to denote by $M$ an abelian group. Let $B$ be an $M$-graded ring and let $I \subseteq B$ a graded ideal. Recall the following graded version of Nakayama. 

\begin{proposition}\label{GradedNakayama}
Suppose that the ideal $I_0$ of $B_0$ is contained in the Jacobson radical of $B_0$.
\begin{assertionlist}
\item\label{GradedNakayama1}
If $E$ is a finite graded $B$-module with $E \otimes_B B/I = 0$, then $E = 0$.
\item\label{GradedNakayama2}
Let $E$ and $F$ be finite graded $B$-modules with $F$ projective as graded $B$-module. Then a map of graded $B$-modules $u\colon E \to F$ is bijective if and only if the induced map $\ubar\colon E \otimes_B B/I \to F \otimes_B B/I$ is bijective.
\end{assertionlist}
\end{proposition}

If $B$ is $\ZZ$-graded with $B_n = 0$ for $n < 0$ and $I = \bigoplus_{n>0}B_n$, then $B/I = B_0$ and the proposition is well known; see, for instance, \cite{Eisenbud_CommutativeAlgebra}. In general, we argue as in \cite[Lemma~3.1.1, Corollary~3.1.2]{Lau_HigherFrames} where the result is proved in a special case of the $\ZZ$-graded case. For the convenience of the reader, we reproduce Lau's proof.

\begin{proof}
We prove \ref{GradedNakayama1}. Let $e_1,\dots,e_r \in E$ be homogeneous generators of $E$ with $e_i$ of degree $m_i \in M$. Since $E = IE$ we can write $e_i = \sum_j a_{ij}e_j$ for $a_{ij} \in I_{m_i - m_j}$. Let $\chi$ be the characteristic polynomial of the matrix $A := (a_{ij})_{i,j}$. As every principal minor of $A$ has degree $0$, we see that all the non-leading coefficients of $\chi$ are in $I_0$. By Cayley-Hamilton, one has $0 = \chi(\id_E) = (1-a)\id_E$ for $a \in I_0$. But $1-a$ is a unit since $I_0$ is in the Jacobson radical of $B_0$. Hence $E = 0$.

Let us show \ref{GradedNakayama2}. If $\ubar$ is bijective, then $u$ is surjective by applying \ref{GradedNakayama1} to $\Coker(u)$. Since $F$ is projective, we find an isomorphism $E \cong F \oplus \Ker(u)$ such that $u$ is the projection onto its first factor. Then $\Ker(u)$ is a finite graded $B$-module with $\Ker(u) \otimes_B B/IB = 0$. Hence $\Ker(u) = 0$ by \ref{GradedNakayama1}, i.e., $u$ is injective.
\end{proof}

We endow $B_0$ with the trivial $M$-grading. Then the maps $B_0 \to B \to B/IB$ of $M$-graded rings, by pullback via
\[
[H_A\backslash (\Spec B/IB)] \to [H_A\backslash (\Spec B)] \to [H_A\backslash (\Spec B_0)],
\]
induce functors
\[
\Vec{[H_A\backslash(\Spec B_0)]} \ltoover{\alpha} \Vec{[H_A\backslash(\Spec B)]} \ltoover{\beta} \Vec{[H_A\backslash(\Spec B/IB)]}.
\]

\begin{theorem}\label{ProjectiveModulesGradedRings}
Suppose that $I_m = B_m$ for all $m \ne 0$ (e.g., if $I$ defines the fixed point locus of $H$ in $\Spec B$, see Remark~\ref{DescribeFixGGm}) and that $(B_0, I_0)$ is a henselian pair. Then $\beta$ is essentially surjective, full, and reflects isomorphisms, and both $\alpha$ and $\beta$ induce bijections on isomorphism classes.
\end{theorem}

\begin{proof}
We have $B_0/I_0 = B/IB$ and both are trivially graded. Since $(B_0, I_0)$ is a henselian pair, $\beta \circ \alpha$ is essentially surjective, full, and reflects isomorphisms (this is \cite[0D4A]{Stacks} in the non-graded case which immediately implies the graded case since $B_0$ and $B/IB$ are both trivially graded). In particular, $\beta \circ \alpha$ induces a bijection on isomorphism classes and $\beta$ is essentially surjective. Let us show that $\beta$ is full and reflects isomorphisms which shows that also $\beta$, and hence $\alpha$, induces a bijection on isomorphism classes.

Let $E$ and $F$ be finite projective graded $B$-modules and let $\ubar\colon E \otimes_B B/IB \lto F \otimes_B B/IB$ be a map of graded $B/IB$-modules. By composition with $E \to E \otimes_B B/IB$ we obtain a map of graded $B$-modules $E \to F \otimes_B B/IB$. Since $E$ is a projective graded $B$-module, this map can be lifted to a map of graded $B$-modules $E \to F$. If $\ubar$ is an isomorphism, then this map is an isomorphism by Proposition~\ref{GradedNakayama}~\ref{GradedNakayama2}.
\end{proof}

%% file: Ch3GBundlesProjectiveLine.tex
\section{Applications of the main theorems}

We now explain two applications of our main results: the classification of $G$-bundles on the projective line (Section~\ref{SecClassGBunPP1}) and the splitting of filtered fiber functors (Section~\ref{SecSplitFiberFunctor}). In both cases we obtain a classification via $G$-bundles on $\B{\GG_m}$. Therefore we start by describing $\Bun_G(\B{\GG_m})$ in Section~\ref{SecGBundleBGGm}.


\subsection{$G$-bundles on the classifying stack of $\GG_m$}\label{SecGBundleBGGm}

Let $S$ be an algebraic space and let $G$ be a flat quasi-compact group scheme over $S$. By Remark~\ref{HomClassifyingStack} the groupoid of $G$-bundles on $\B{\GG_{m,S}}$ is the groupoid of pairs $(P, \mu)$, where $P$ is a $G$-bundle on $S$ and $\mu\colon \GG_{m,S} \to \Autline_G(P)$ is a cocharacter of the strong inner form of $G$ attached to $P$. A morphism $(\Escr,\mu) \to (\Escr',\mu')$ is a morphism of $G$-bundles intertwining $\mu$ and $\mu'$.

If $G$ is affine and smooth over $S$, \eqref{EqDescribeGBundleClassStack2} and Proposition~\ref{BunGClassStackAlgebraic} show that
\[
\Bunline_G(\B{\GG_{m,S}}) = [G\backslash \Homline_{\textup{$S$-Grp}}(\GG_{m,S},G)]
\]
is an algebraic stack locally of finite presentation over $S$.

We also have the following Tannakian description of $\Bun_G(\B{\GG_m})$. 

\begin{remark}\label{BunGBGGmTannaka}
Suppose that $G$ is obtained by base change $S \to S_0$ from a flat and quasi-affine group scheme $G_0$ over $S_0$ such that $\B{G_0}$ satisfies the resolution property (see Example~\ref{BGResolution}). Then $\Bun_G(\B{\GG_m}) = \Hom_{S_0}(\B{\GG_{m,S}},\B{G_0})$ is equivalent to the groupoid of exact $\Oscr_{S_0}$-linear strong monoidal functors $\Rep(G_0) \to \Rep(\GG_{m,S})$ (Proposition~\ref{TannakaGBundles}). Since $\Rep(\GG_{m,S})$ can be identified with the category of $\ZZ$-graded vector bundles over $S$, we see that $\Bun_G(\B{\GG_m})$ is equivalent to the groupoid of graded fiber functors of $\Rep(G_0)$ over $S$ in the sense of \cite[IV, \S1]{SR_Tannaka} (see also \cite[Section~3.2]{Ziegler_FFF}).
\end{remark}

%
%
%
%
%
%


\subsection{Classification of $G$-bundles over the projective line}\label{SecClassGBunPP1}

As an application we obtain a new proof of the classification of $G$-bundles over $\PP^1_k$ in terms of $G$-bundles over $\B{\GG_m}$:

\begin{theorem}\label{ClassifyGBundlesPP1}
Let $k$ be a field and let $G$ be a smooth affine group scheme over $k$ whose identity component $G^0$ is reductive. Then one has a bijection
\begin{equation}\label{EqBijGBunPP1}
H^1(\B{\GG_{m,k}},G) \cong H^1(\PP^1_k,G).
\end{equation}
\end{theorem}

\begin{proof}
As $G$ is defined over a field, $G/G^0$ is automatically finite. By Remark~\ref{RemGeomRedGrp}~\ref{RemGeomRedGrp1} we see that $G$ is geometrically reductive.
Consider $X = \AA^2_k$ with its standard $\GG_m$-action. Then one has bijective maps
\[\xymatrix{
H^1(\B{\GG_{m,k}},G) & H^1([\GG_m\backslash \AA^2],G) \ar[l]_-{\sim} \ar[r]^-{\sim} & H^1(\PP^1_k,G).
}\]
Here the first map is a special case of \eqref{EqH1TrivialHensel} and one has the second bijection by Proposition~\ref{ExtendGBundleOpenImmersion}~\ref{ExtendGBundleOpenImmersion2} using that $\PP^1_k = [(\AA^2_k \setminus \{0\})/\GG_m]$.
\end{proof}

\begin{remark}\label{ExplicitDescriptionGBunPP1}
One can make \eqref{EqBijGBunPP1} more explicit as follows. As recalled above, a $G$-bundle over $\B{\GG_{m,k}}$ is given by a $G$-bundle $P$ and by a map of group schemes $\mu\colon \GG_{m,k} \to \Autline_G(P)$.

Let $f\colon \AA_k^2 \setminus \{0\} \to \Spec k$ be the structure morphism. By Example~\ref{TrivialHenselian}, the corresponding $G$-bundle on $[\GG_{m,k}\backslash (\AA_k^2 \setminus \{0\})]$ is given by the $G$-bundle $\Ptilde := f^*P$ over $\AA_k^2 \setminus \{0\}$ and the $\GG_m$-equivariant structure defined by $\mgtilde\colon \GG_{m,k} \to \Autline_G(\Ptilde)$ obtained by base change from $\mu$. Since the geometric line bundle $\AA^2_k \setminus \{0\} \to \PP^1_k$ is $\Oscr(-1)$, the $G$-bundle over $\PP^1_k$ corresponding to $(P,\mu)$ is given by the image of the isomorphism class of $\Oscr(-1)$ under the composition
\[
H^1(\PP^1_k,\GG_{m,k}) \ltoover{\mgtilde_*} H^1(\PP^1_k,\Autline_G(P)) \vartoover{35}{(\Theta^P)^{-1}} H^1(\PP^1_k,G),
\] 
where $\Theta^P$ is the twist by $P$.
\end{remark}

Since $G$ is smooth over $k$, every $G$-bundle is \'etale locally trivial. Moreover we can also characterize all $G$-bundles on $\PP^1_k$ that are Zariski locally trivial:

\begin{proposition}\label{ZarLocallyTrivPP1}
Let $G$ be as in Theorem~\ref{ClassifyGBundlesPP1}, let $E$ be a $G$-bundle over $\PP^1_k$, and let $(P,\mu)$ be the corresponding $G$-bundle over $\B{\GG_{m,k}}$. Let $x \in \PP^1(k)$ be a $k$-rational point. Then the following assertions are equivalent.
\begin{equivlist}
\item\label{ZarLocallyTrivPP1i}
The $G$-bundle $E$ is Zariski locally trivial.
\item\label{ZarLocallyTrivPP1ii}
The $G$-bundle $P$ on $k$ is trivial.
\item\label{ZarLocallyTrivPP1iii}
The pullback $x^*E$ is trivial.
\end{equivlist}
\end{proposition}

\begin{proof}
We use the notation introduced in Remark~\ref{ExplicitDescriptionGBunPP1}. We first claim that we have $x^*E \cong P$. Indeed $x$ can be lifted to a $k$-rational point $\xtilde$ of $\AA^2 \setminus \{0\}$. Using the notation of Remark~\ref{ExplicitDescriptionGBunPP1}, we find $x^*E = \xtilde^*\Ptilde = P$.

Clearly, \ref{ZarLocallyTrivPP1i} implies \ref{ZarLocallyTrivPP1iii}. The above claim shows that \ref{ZarLocallyTrivPP1iii} implies \ref{ZarLocallyTrivPP1ii}. 

It remains to show that \ref{ZarLocallyTrivPP1ii} implies \ref{ZarLocallyTrivPP1i}. If $P$ is trivial, then Remark~\ref{ExplicitDescriptionGBunPP1} shows that $E = \mgtilde_*(\Oscr(-1))$. Since $\Oscr(-1)$ is Zariski locally trivial, $E$ is Zariski locally trivial.
\end{proof}

As a special case one also obtains the following more classical version (a slight generalization of \cite{BiswasNagaraj}) of classification of Zariski locally trivial $G$-bundles over $\PP^1_k$ whose formulation avoids $G$-bundles over an algebraic stack.

\begin{corollary}\label{ClassifyGBundlesPP1Classical}
Let $k$ be a field and let $G$ be a smooth affine group scheme over $k$ whose identity component is reductive. Let $P$ be a $G$-bundle over $\PP^1_k$ that is trivial over some $k$-rational point of $\PP^1_k$. Then there exists a cocharacter $\lambda\colon \GG_{m,k} \to G$ (defined over $k$) such that the isomorphism class of $P$ is obtained as the image of the $\GG_m$-bundle $\Oscr(-1)$ under the map $H^1(\PP^1_k,G) \to H^1(\PP^1_k,G)$ induced by $\lambda$.
\end{corollary}


\begin{remark}\label{Pi1PP1}
One can also apply Theorem~\ref{ClassifyGBundlesPP1} to finite \'etale group schemes $G$ over $k$. Then every cocharacter $\GG_m$ of any strong inner form of $G$ is trivial and pullback along $\Spec k \to \B{\GG_{m,k}}$ yields an equivalence $\Bun_G(\B{\GG_{m,k}}) \iso \Bun_G(k)$. Therefore Theorem~\ref{ClassifyGBundlesPP1} yields in this case a bijection
\[
H^1(k,G) \cong H^1(\PP^1_k,G).
\]
Hence every connected finite \'etale Galois covering of $\PP^1_k$ arises by base change from a finite Galois extension of $k$ which implies the well known isomorphism
\[
\pi_1(\PP^1_k) \iso \Gal(k^{\sep}/k)
\]
of \'etale fundamental groups.
\end{remark}


\subsection{Splitting of filtered fiber functors}\label{SecSplitFiberFunctor}

Let $S_0$ be a qcqs scheme and let $G_0 \to S_0$ be a flat quasi-affine group scheme such that $\B{G_0}$ satisfies the resolution property (see Example~\ref{BGResolution} for examples). Let $S$ be an $S_0$-scheme and set $G := S \times_{S_0} G_0$. We now consider $G$-bundles on $\Xcal := [\GG_{m,S}\backslash \AA_S^1]$. By Proposition~\ref{TannakaGBundles}, sending a morphism $f$ to the pullback functor $f^*$ yields an equivalence of the groupoid
\[
\Bun_G(\Xcal) = \Hom_S(\Xcal,\B{G}) = \Hom_{S_0}(\Xcal,\B{G_0})
\]
and the groupoid of exact $\Oscr_{S_0}$-linear strong monoidal functors $\Rep(G_0) \to \Vec{\Xcal}$. Since $\Vec{\Xcal}$ can be identified with the category of vector bundles on $S$ endowed with a separated and exhaustive increasing $\ZZ$-filtration by subbundles (Example~\ref{VectorbundlesoverAA1ModGGm}), it follows that $\Bun_G(\Xcal)$ can be identified with the groupoid of filtered fiber functors of $\Rep(G_0)$ over $S$.

The zero section of $S \to \AA^1_S$ yields a closed immersion
\[
\B{\GG_{m,S}} = [\GG_{m,S}\backslash S] \ltoover{i} [\GG_{m,S}\backslash \AA^1_S].
\]
Pullback of $G$-bundles via $i$ defines a functor
\begin{equation}\label{EqSplittingFiltered}
i^*\colon \Bun_G([\GG_{m,S}\backslash \AA^1_S]) \lto \Bun_G(\B{\GG_{m,S}}).
\end{equation}
It attaches to every filtered fiber functor of $\Rep(G_0)$ over $S$ its associated graded fiber functor (Remark~\ref{BunGBGGmTannaka}) induced by the exact strong monoidal functor that sends a filtered vector bundle $(\Escr, (\Fil^i\Escr)_{i \in \ZZ})$ over $S$ to the $\ZZ$-graded vector bundle $\bigoplus_i \Fil^i\Escr/\Fil^{i-1}\Escr$.

The structure map $\AA^1_S \to S$ induces a map $p\colon [\GG_{m,S}\backslash \AA^1_S] \to \B{\GG_{m,S}}$ and hence via pullback an exact strong monoidal functor
\begin{equation}\label{EqFilteredGraded}
p^*\colon \Bun_G(\B{G}) \lto \Bun_G([\GG_{m,S}\backslash \AA^1_S]).
\end{equation}
induced by the exact strong monoidal functor that sends a graded vector bundle $(\Fscr = \bigoplus_i \Fscr_i)$ over $S$ to the filtered vector bundle $(\Fscr, \Fil^i\Fscr := \bigoplus_{j\leq i}\Fscr_j)$. Clearly one has that $i^* \circ p^* \cong \id$.

A filtered fiber functor of $\Rep(G_0)$ over $S$ is then called \emph{splittable}, if it lies in the essential image of $p^*$.

We obtain as another direct application of Theorem~\ref{LiftHenselian} the following result which for smooth affine group schemes over fields is a special case of \cite[Theorem~4.16]{Ziegler_FFF}.

\begin{corollary}\label{SplittingFFF}
Let $S_0$, $G_0$, $S$, and $G$ as above. Suppose that $S$ is affine and that $G$ is smooth and quasi-affine over $S$. Then every filtered fiber functor of $\Rep(G_0)$ over $S$ is splittable by a graded fiber functor which is uniquely determined up to isomorphism.
\end{corollary}

\begin{proof}
By Theorem~\ref{LiftHenselian}, the functor \eqref{EqSplittingFiltered} is essentially surjective and full. This implies that \eqref{EqFilteredGraded} induces a bijection $H^1(\B{\GG_{m,S}},G) \iso H^1([\GG_{m,S}\backslash \AA^1_S],G)$.
\end{proof}

Note that the bijection $H^1([\GG_{m,S}\backslash \AA^1_S],G) \iso H^1(\B{\GG_{m,S}},G)$ holds for arbitrary smooth quasi-affine group schemes $G$ over an affine scheme $S$. The resolution property is needed only to identify both sides with filtered, respectively graded, fiber functors.

\begin{remark}\label{SmoothnessNecessary}
Corollary~\ref{SplittingFFF} also shows that the smoothness hypothesis on $G$ is necessary for the uniqueness part in Theorem~\ref{LiftHenselian} as Saavedra Rivano has given in \cite[IV.2.2.3.1]{SR_Tannaka} examples of a finite type non-smooth group scheme $G$ over a non-perfect field $k$ and a filtered fiber functor of $\Rep(G)$ that is not splittable over $k$ (but only over a non-separable extension of $k$).
\end{remark}

%


%
